\documentclass[a4paper, 11pt]{article}
\usepackage[english]{babel}
\usepackage{amsmath,amssymb,graphicx,stmaryrd,enumerate,mathrsfs,bbm,xcolor,theorem,a4wide}

\usepackage[hidelinks=true]{hyperref}

\newcommand{\assign}{:=}
\newcommand{\longdownarrow}{{\mbox{\rotatebox[origin=c]{-90}{$\longrightarrow$}}}}
\newcommand{\longuparrow}{{\mbox{\rotatebox[origin=c]{90}{$\longrightarrow$}}}}
\newcommand{\mathd}{\mathrm{d}}

\newcommand{\tmop}[1]{\ensuremath{\operatorname{#1}}}
\newcommand{\tmtextit}[1]{{\itshape{#1}}}
\newtheorem{theorem}{Theorem}[section]

\newtheorem{lemma}[theorem]{Lemma}
\newtheorem{remark}[theorem]{Remark}

\newcommand{\zzone}{\text{\resizebox{.7em}{!}{\includegraphics{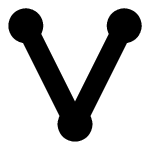}}}}
\newcommand{\zztwo}{\text{\resizebox{.7em}{!}{\includegraphics{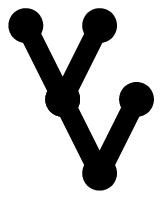}}}}
\newcommand{\zzthreereso}{\text{\resizebox{.7em}{!}{\includegraphics{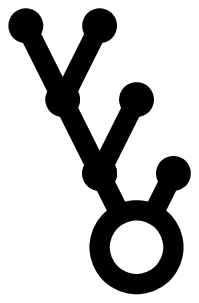}}}}
\newcommand{\zzfour}{\text{\resizebox{1em}{!}{\includegraphics{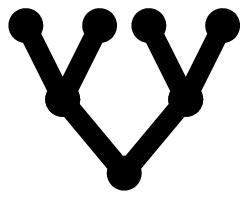}}}}
\newcommand{\CC}{\mathcal{C}}

\newcommand{\LL}{\mathcal{L}}
\newcommand{\dd}{\mathrm{d}}

\newcommand{\para}{\mathord{\prec}}
\newcommand{\lpara}{\mathord{\succ}}
\newcommand{\mpara}{\mathord{\prec\!\!\!\prec}}
\newcommand{\reso}{\mathord{\circ}}

\makeindex             % used for the subject index
\begin{document}

\title{An introduction to singular SPDEs}
% Use \titlerunning{Short Title} for an abbreviated version of
% your contribution title if the original one is too long
\author{
  Massimiliano Gubinelli\\
  Hausdorff Center for Mathematics\\
   \& Institute for Applied Mathematics\\
   Universit{\"a}t Bonn \\
  \texttt{gubinelli@iam.uni-bonn.de}
  \and
  Nicolas Perkowski\thanks{Financial support by the DFG via Research Unit FOR 2402 is gratefully acknowledged.} \\
  Institut f\"ur Mathematik \\
  Humboldt--Universit\"at zu Berlin \\
  \texttt{perkowsk@math.hu-berlin.de}
}

%
% Use the package "url.sty" to avoid
% problems with special characters
% used in your e-mail or web address
%
\maketitle

%\abstract*{Each chapter should be preceded by an abstract (10--15 lines long) that summarizes the content. The abstract will appear \textit{online} at \url{www.SpringerLink.com} and be available with unrestricted access. This allows unregistered users to read the abstract as a teaser for the complete chapter. As a general rule the abstracts will not appear in the printed version of your book unless it is the style of your particular book or that of the series to which your book belongs.
%Please use the 'starred' version of the new Springer \texttt{abstract} command for typesetting the text of the online abstracts (cf. source file of this chapter template \texttt{abstract}) and include them with the source files of your manuscript. Use the plain \texttt{abstract} command if the abstract is also to appear in the printed version of the book.}

\begin{center}
\emph{Dedicated to Michael R{\"o}ckner on the occasion of his 60th
  Birthday.}
\end{center}

\abstract{We review recent results on the analysis of singular stochastic partial
  differential equations in the language of paracontrolled distributions.}

\section{Introduction}

In recent years there has been much progress in the mathematical understanding
of certain non-linear random PDEs which are not well posed in the classical
analytic or probabilistic theory but which become amenable to rigorous
analysis as soon as specific non-linear properties of the randomness are taken
into account. In the related literature it has become common to refer to such
equations as \tmtextit{singular SPDEs}, mainly to distiguish them from
standard SPDEs. The difference is that singular SPDEs can be posed only in
small subspaces of the standard function spaces (e.g. H{\"o}lder, Sobolev or
even Besov spaces) and that the operations involved in such equations
sometimes require \tmtextit{renormalisation}. Renormalisation in this context
can be understood as the fact that only specific non-linearities can be formed
meaningfully and that, in order to do so, subtractions of infinite quantities
are often needed.

The aim of this short review is to present some of the ideas underlying
singular SPDEs and their pathwise analysis and to make the connection with
other parts of probability theory and mathematical physics: to renormalisation
group theory and to scaling limits of interacting particle models. An
important motivation to consider singular SPDEs is indeed that in some cases
they appear in the description of universal large scale
fluctuations for certain spatially extended probabilistic models. Universality
here means that irrespective of most specific features of the model its large
scale fluctuations are described by a generic equation that usually can be
guessed by first principles and then hopefully confirmed by rigorous analysis.
One of the main open problem in singular SPDEs is to enlarge the spectrum of
models for which the universality can be rigorously proven.

%Currently we dispose of three main approaches in order to study singular
%equations: \tmtextit{regularity structures}, \tmtextit{paracontrolled
%distributions} and \tmtextit{the renormalisation group approach}. Regularity
%structures have been introduced by M. Hairer in his remarkable
%work~{\cite{hairer_theory_2014}} where they were in particular used to give, for the first time, a solution theory for the dynamic $\Phi^4_3$ model.
Currently we dispose of four main approaches in order to study singular
equations: \tmtextit{regularity structures}, \tmtextit{paracontrolled
distributions}, \tmtextit{the renormalisation group approach}, and Otto's and Weber's rough path based approach. Apart from the renormalisation group approach, the other three techniques are all
inspired by T.~Lyons's \emph{rough path}
theory~{\cite{lyons_differential_1998,lyons_system_2002,lyons_differential_2007}} and by the related notion of \emph{controlled paths}~{\cite{gubinelli_controlling_2004,friz_course_2014}}. 
 Regularity
structures have been introduced by M. Hairer in his remarkable
work~{\cite{hairer_theory_2014}} where they were in particular used to give, for the first time, a solution theory for the dynamic $\Phi^4_3$ model. Regularity
structures allow a detailed description of the \tmtextit{local} action of
distributions on test functions in terms of a given \tmtextit{model} which
usually is constructed from certain non-linear features of an underlying
random process. Paracontrolled distributions have been introduced by the
authors together with P.~Imkeller~{\cite{gubinelli_paracontrolled_2015}}, more
or less at the same time as M.~Hairer was developing his theory, as a tool to
describe the ``spectral'' features of a function (or distribution) in terms of
simpler objects, very similar to Hairer's models. Some time after,
Kupiainen~{\cite{kupiainen_renormalization_2016-1}} observed that the
renormalisation group approach can be also used to analyse singular equations
by decomposing random fields in a multiscale fashion and by introducing
recursive equations for each scale. The most recent approach is due to Otto and Weber~\cite{otto_quasilinear_2016}. In their approach a semigroup is used to provide a multiscale resolution of various singular objects and the scale parameter is handled in the spirit of the time parameter in rough path theory. In this review we will not address the the connections of
paracontrolled distributions with the the other techniques.

We will illustrate the analysis of singular SPDEs on a series of models:
\begin{enumerate}
  \item The 1d generalised Stochastic Burgers equation (gSBE)
  \begin{equation}
    \partial_t u (t, x) = \Delta u (t, x) + G (u (t, x)) \partial_x u (t, x) +
    \xi (t, x), \qquad t \geqslant 0, x \in \mathbb{T}, \label{eq:gSBE}
  \end{equation}
  where $u\colon \mathbb{R}_+ \times \mathbb{T} \rightarrow \mathbb{R}$ and
  $G : \mathbb{R} \rightarrow \mathbb{R}$ is a 
  smooth function.
  
  \item The 1d Kardar--Parisi--Zhang equation (KPZ)
  \begin{equation}
    \partial_t h (t, x) = \Delta h (t, x) + (\partial_x h (t, x))^2 - C + \xi
    (t, x), \qquad t \geqslant 0, x \in \mathbb{T} \label{eq:KPZ}
  \end{equation}
  with $h\colon \mathbb{R}_+ \times \mathbb{T} \rightarrow \mathbb{R}$ and $C
  \in \mathbb{R}$ and the related conservative stochastic Burgers equation
  (CSBE) 
  \begin{equation}
    \partial_t u (t, x) = \Delta u (t, x) + \partial_x (u (t, x)^2) +
    \partial_x \xi (t, x), \qquad t \geqslant 0, x \in \mathbb{T},
    \label{eq:CSBE}
  \end{equation}
  where $u = \partial_x h\colon \mathbb{R}_+ \times \mathbb{T} \rightarrow
  \mathbb{R}$.
  
  \item The dynamic $\Phi^4_d$ model or stochastic quantisation equation ($d =
  2, 3$) (SQE)
  \begin{equation}
    \partial_t \varphi (t, x) = \Delta \varphi (t, x) - \varphi (t, x)^3 + C \varphi (t, x) + \xi (t, x), \qquad t \geqslant 0,
    x \in \mathbb{T}^d, \label{eq:SQE}
  \end{equation}
  where $\varphi\colon \mathbb{R}_+ \times \mathbb{T}^d \rightarrow \mathbb{R}$
  and $C \in \mathbb{R}$.
  
  \item The generalised two--dimensional parabolic Anderson model (gPAM)
  \begin{equation}
    \partial_t u (t, x) = \Delta u (t, x) + G (u (t, x)) \xi (x) - C G' (u (t,
    x)) G (u (t, x)), \qquad t \geqslant 0, x \in \mathbb{T}^2,
    \label{eq:gPAM}
  \end{equation}
  where $u\colon \mathbb{R}_+ \times \mathbb{T}^2 \rightarrow \mathbb{R}$, $G :
  \mathbb{R} \rightarrow \mathbb{R}$ is a smooth function, and $C \in
  \mathbb{R}$.
\end{enumerate}
In all these examples $\mathbb{T}^d =\mathbb{R}^d /\mathbb{Z}^d$ is the
$d$-dimensional torus and $\xi$ denotes a Gaussian white noise (space--time or
space dependent only, according to the model). The specific choice of the
dimensionality and/or space--time dependence of the noise is related to the
degree of singularity of the equation and is motivated by the fact that in
this review we will use the specific language of paracontrolled distributions
which has the advantage of requiring very little background and is more
directly related to standard PDE theory than the other approaches.

This choice has the drawback that we will not address the discussion of
natural models like the generalised form of the KPZ equation, given by
\[ \partial_t h (t, x) = \Delta h (t, x) + G (h (t, x)) (\partial_x h (t,
   x))^2 + F (h (t, x)) \xi (t, x), \qquad t \geqslant 0, x \in \mathbb{T},
\]
which is within reach of regularity structure theory but still out of reach
for paracontrolled distributions. The generalised KPZ equation is in a way the
``ultimate'' singular SPDE and its singularities are very challenging. Its
analysis via regularity structures requires a great deal of work and a deeper
understanding of the algebraic and analytic structures underlying the
construction of an appropriate model and of the
renormalisation~{\cite{bruned_algebraic_2016,chandra_analytic_2016}}. 

As we already remarked the main difficulty shared by all singular SPDEs is the
presence of non-linear operations which are not well defined in classical
function spaces. This difficulty blocks the analysis from the very beginning
because it is not even possible to \tmtextit{formulate} the equation
rigorously. All the equations we wrote in the introduction are classically
ill-posed and the notation is just used in an informal and suggestive way.
Indeed most of them have to be modified to take into account
\tmtextit{renormalisation} which means that the constant $C$ that appear in
the equations is actually to be read as an infinite quantity and not a real
number. A standard approach to make the analysis rigorous is to construct a
series of approximate problems whose solutions converge to a well-defined
limit. The characterisation of this limit and its independence of the details
of the approximation procedure will constitute a rigorous definition of the
solution to a singular SPDE. Taking as an example the 2d gPAM, we can
(loosely) state a typical convergence result:

\begin{theorem}
  Let $G : \mathbb{R} \rightarrow \mathbb{R}$ be a smooth function,
  $\varepsilon > 0$, and let $u_{\varepsilon} : \mathbb{R}_+ \times
  \mathbb{T}^2 \rightarrow \mathbb{R}$ be the unique classical solution of
  the Cauchy problem
  \[ \partial_t u_{\varepsilon} (t, x) = \Delta u_{\varepsilon} (t, x) + G
     (u_{\varepsilon} (t, x)) \xi_{\varepsilon} (x) - G_2 (u_{\varepsilon} (t,
     x)) c_{\varepsilon}, \qquad t \geqslant 0, x \in \mathbb{T}^2, \]
  where $\xi_{\varepsilon}$ is a smooth approximation of the space white noise
  $\xi : \mathbb{T}^2 \rightarrow \mathbb{R}$ by convolution,
  $c_{\varepsilon} \in \mathbb{R}$, and $G_2 (u) = G' (u) G (u)$. Then there
  exists a choice for $c_{\varepsilon}$ such that locally in time the sequence
  $u_{\varepsilon}$ converges, as $\varepsilon \rightarrow 0$, to a function
  $u$ that does not depend on the exact form of the convolutional
  approximation of $\xi_{\varepsilon}$.
\end{theorem}

This formulation is quite minimal, in reality the analysis which has to be put
forward to prove this kind of results give as a byproduct also much more
detailed information about the function $u$. In particular it is possible to
characterise $u$ via a standard PDE for a different unknown and the
assumptions on the approximations of the noise $\xi$ become assumptions about
the convergence of certain non-linear functionals of $\xi_{\varepsilon}$. The
general scheme underlying the convergence is the following:
\[ \begin{array}{lllllll}
     \xi_{\varepsilon} & \overset{J}{\longleftrightarrow} & \Xi_{\varepsilon}
     & \overset{\Phi}{\longmapsto} & U_{\varepsilon} & \overset{\Pi}{\longmapsto} &
     u_{\varepsilon}\\
     \longdownarrow &  & \longdownarrow &  & \longdownarrow &  &
     \longdownarrow\\
     \xi & \longmapsfrom & \Xi & \overset{\Phi}{\longmapsto} & U & \overset{\Pi}{\longmapsto}
     & u
   \end{array} \]
Here the vertical arrows represent limits for $\varepsilon \rightarrow 0$ (in
appropriate function spaces). The upper row features $\Xi_{\varepsilon} = J
(\xi_{\varepsilon})$, an injective collection of non-linear quantities
constructed from the approximate data $\xi_{\varepsilon}$ (one could also
consider the initial condition to be part of this data but we will refrain
from doing so), and $U_{\varepsilon} = \Phi (\Xi_{\varepsilon})$, an
\tmtextit{enhanced} notion of solution from which one can recover the
classical solution $u_{\varepsilon}$ through a continuous projection $\Pi$. The
bottom row describes the situation after the limit $\varepsilon \rightarrow 0$
has been taken. The \tmtextit{enhanced data} $\Xi$ still determines the limit
noise $\xi$, however the reverse in not true and different approximation
procedures for the same $\xi$ can lead to different values of $\Xi$. But the
remaining relations are preserved: from $\Xi$ we can still recover an enhanced
notion of solution $U$ through the same continuous \tmtextit{solution map}
$\Phi : \Xi \mapsto U$, and from $U$ we obtain $u$ by a projection in such a
way that the convergence $u_{\varepsilon} \rightarrow u$ holds by the
continuity of the solution map and of the projection. In the limit the
situation is more complex than before passing to the limit. As mentioned
before a different approximation procedure $\tilde{\xi}_{\varepsilon}
\rightarrow \xi$ can lead to different enhanced data $\tilde{\Xi}$ in the
limit, and as a consequence of the continuity of $\Phi$ to a different
limiting solution $\tilde{u}$. This situation is depicted in the graph below.
\[ \begin{array}{lllllll}
     \xi_{\varepsilon} & \overset{J}{\longleftrightarrow} & \Xi_{\varepsilon}
     & \overset{\Phi}{\longmapsto} & U_{\varepsilon} & \overset{\Pi}{\longmapsto} &
     u_{\varepsilon}\\
     \longdownarrow &  & \longdownarrow &  & \longdownarrow &  &
     \longdownarrow\\
     \xi & \longmapsfrom & \Xi & \overset{\Phi}{\longmapsto} & U & \overset{\Pi}{\longmapsto}
     & u\\
     = &  & \neq &  & \neq &  & \neq\\
     \xi & \longmapsfrom & \tilde{\Xi} & \overset{\Phi}{\longmapsto} &
     \tilde{U} & \overset{\Pi}{\longmapsto} & \tilde{u}\\
     \longuparrow &  &  &  & \longuparrow &  & \longuparrow\\
     \tilde{\xi}_{\varepsilon} & \overset{J}{\longleftrightarrow} &
     \tilde{\Xi}_{\varepsilon} & \overset{\Phi}{\longmapsto} &
     \tilde{U}_{\varepsilon} & \overset{\Pi}{\longmapsto} & \tilde{u}_{\varepsilon}
   \end{array} \]
So the relation between the data $\xi$ and the solution is not well defined
unless some information on the non-linear features $\Xi$ is provided as
additional input.

Proving a convergence result therefore requires two conceptually different
steps:
\begin{enumerate}
  \item define and analyse the continuity of the map $\Phi$, where an
  important task is to identify its domain and co--domain;
  
  \item prove the convergence of the enhanced data $\Xi_{\varepsilon}
  \rightarrow \Xi$.
\end{enumerate}
The convergence of the enhanced data is, in most of the applications, a purely
probabilistic step which sometimes requires development of efficient tools but
for which the main tools are classical and already present in the
probabilistic literature for a long time. Some keywords in this context are
Gaussian analysis, hypercontractivity, Besov embeddings, almost sure
regularity of stochastic processes, martingale theory, Wick products, and
chaos expansions. We would like to concentrate our discussion to the analytic
part of the theory involving the construction of the enhanced spaces which
constitute the domain and co-domain of the solution map $\Phi$.

In Sections~\ref{sec:paraproducts} and~\ref{sec:paracontrolled} we present the basic ideas and analytic ingredients for paracontrolled distributions. Section~\ref{sec:renormalization} briefly discusses the need for renormalisation of the enhanced data and how this renormalisation translates in the equation. In Section~\ref{sec:higherorder} we sketch recent work of Bailleul and Bernicot on higher order expansions via paraproducts. Section~\ref{sec:universality} is dedicated to convergence results for singular SPDEs and we illustrate how to derive the Hairer-Quastel weak universality principle for the KPZ equation using paracontrolled distributions. In Sections~\ref{sec:operators1},~\ref{sec:operators2} we will see that paracontrolled distributions can not only be used to study singular SPDEs, but as noted by Cannizzaro and Chouk respectively Allez and Chouk they also allow us to construct certain singular operators.

Finally let us point out that there are many fascinating recent results that are based on paracontrolled distributions and that we have to omit here due to space constraints. To name just a few: non-explosion results for the dynamic $\Phi^4_3$ model \cite{mourrat_global_2016}, the KPZ equation~\cite{gubinelli_kpz_2017}, and the multi-component KPZ equation \cite{funaki_coupled_2016}, a formulation of paracontrolled distributions that allows to study equations on manifolds \cite{bailleul_heat_2016}, convergence results for discrete dynamics \cite{zhu_approximating_2014, zhu_wong_2015, gubinelli_kpz_2017, chouk_invariance_2016}, a solution theory for quasilinear equations \cite{bailleul_quasilinear_2016}, nonlinear extensions of paraproducts \cite{furlan_paracontrolled_2016}, and a support theorem for gPAM \cite{chouk_support_2014} -- not to mention all the exciting results that have been shown in the setting of regularity structures or the other approaches.

\section{Paraproducts}\label{sec:paraproducts}

In order to develop the paracontrolled analysis of the solution map we will
introduce in this section Besov spaces and paraproducts.
See~{\cite{bahouri_fourier_2011}} for details. Let $\mathcal{S}'
(\mathbb{T}^d)$ denote the Schwartz space of distributions on $\mathbb{T}^d
=\mathbb{R}^d /\mathbb{Z}^d$. Then any element $f \in \mathcal{S}'
(\mathbb{T}^d)$ can be decomposed as
\[ f = \sum_{i \geqslant - 1} \Delta_i f, \]
where the sum is over $i = - 1, 0, 1, 2, \ldots$ and $\Delta_i f$ are smooth
functions whose Fourier support is contained in a ball $\mathcal{B} \subseteq
\mathbb{R}^d$ for $i = - 1$ and in rescaled dyadic annuli $2^i \mathcal{A}$
for $i \geqslant 0$. The operators $\Delta_i : f \mapsto \Delta_i f$ are
called Littlewood--Paley operators and can be constructed to enjoy nice
analytic properties, for example they satisfy the Bernstein inequalities
\[ \| D^{\alpha} \Delta_i f \|_{L^p (\mathbb{T}^d)} \lesssim 2^{i | \alpha |
   + i \left( \frac{d}{q} - \frac{d}{p} \right)} \| \Delta_i f \|_{L^q
   (\mathbb{T}^d)}, \qquad i \geqslant - 1, \qquad p \geqslant q, \]
where $\alpha$ is a $d$-dimensional multiindex and $D^{\alpha}$ denotes the
related mixed derivative of order $| \alpha |$. The (inhomogeneous) Besov
space $B^{\alpha}_{p, q}$ is defined as the set of all distributions $f \in
\mathcal{S}' (\mathbb{T}^d)$ such that the norm
\[ \| f \|_{B^{\alpha}_{p, q}} = \| (\| 2^{i \alpha} \Delta_i f \|_{L^p
   (\mathbb{T}^d)})_{i \geqslant - 1} \|_{\ell^q (\mathbb{Z})} \]
is finite. In the analysis of singular SPDEs we will mainly use the
H{\"o}lder--Besov spaces $\mathcal{C}^{\alpha} = B^{\alpha}_{\infty, \infty}$
to get rid of the integrability exponents. This turns out to be very
convenient for non-linear estimates because $\| f g \|_{L^{\infty}} \leqslant
\| f \|_{L^{\infty}} \| g \|_{L^{\infty}}$ and similarly for $\ell^{\infty}$
which is not true for the $L^p$ respectively $\ell^p$ norms with $p < \infty$.
To gain an intuitive understanding of the $\mathcal{C}^{\alpha}$ spaces it is
useful to note that for $\alpha \in \mathbb{R}_+ \setminus \mathbb{Z}$ the
space $\mathcal{C}^{\alpha}$ consists exactly of the $\lfloor \alpha \rfloor$
times continuously differentiable functions from $\mathbb{T}^d \rightarrow
\mathbb{R}$ for which the classical increment-based $(\alpha - \lfloor \alpha
\rfloor)$-H{\"o}lder norm of all partial derivatives of order $\lfloor \alpha
\rfloor$ is finite. And roughly speaking for $\alpha \in (- 1, 0)$ every $f
\in \mathcal{C}^{\alpha}$ is the distributional derivative of some $F \in
\mathcal{C}^{\alpha + 1}$, and similarly distributions of lower regularity are
higher order derivatives of H{\"o}lder continuous functions.

The Littlewood--Paley decomposition induces a natural decomposition of
products of Besov functions in terms of paraproducts. Given $f, g$ we have
\[ fg = \sum_{i, j \geqslant - 1} \Delta_i f \Delta_j g = f \para g + f \reso
   g + f \lpara g, \]
where the \tmtextit{paraproducts} $f \para g$ and $f \lpara g$ and the
\tmtextit{resonant product} $f \reso g$ by are defined, respectively, by
\[ f \para g = g \lpara f \assign \sum_{i \geqslant - 1} (\Delta_{< i - 1} f)
   \Delta_i g, \qquad f \reso g \assign \sum_{i, j : | i - j | \leqslant 1}
   \Delta_i f \Delta_j g, \]
and where we introduced the notation $\Delta_{< k} f = \sum_{\ell < k}
\Delta_{\ell} f$. Paraproducts are continuous bilinear operations on the
following function spaces
\[ \ast \para \ast : \mathcal{C}^{\alpha} \times \mathcal{C}^{\beta}
   \rightarrow \mathcal{C}^{\beta + \alpha}, \qquad \qquad \alpha \leqslant 0,
   \beta \in \mathbb{R}, \]
\[ \ast \para \ast : L^{\infty} \times \mathcal{C}^{\beta} \rightarrow
   \mathcal{C}^{\beta}, \qquad \beta \in \mathbb{R}, \]
while the resonant product is well defined only if $\alpha + \beta > 0$ and in
that case it is a continuous bilinear operator
\[ \ast \reso \ast : \mathcal{C}^{\alpha} \times \mathcal{C}^{\beta}
   \rightarrow \mathcal{C}^{\beta + \alpha} . \]
In particular we see that the usual product  can
be extended by continuity from smooth functions to H{\"o}lder--Besov
distributions as a bilinear map
\[ (f, g) \mapsto fg : \mathcal{C}^{\alpha} \times \mathcal{C}^{\beta}
   \rightarrow \mathcal{C}^{\min (\alpha, \beta)} \]
provided $\alpha + \beta > 0$.

Singular SPDEs are characterised by this condition not being satisfied in the
non-linear terms. The products then become problematic and undefined for
general inputs. To explain the difficulty let us note that given $\alpha,
\beta \in \mathbb{R}$ with $\alpha + \beta < 0$ and given smooth functions
$(f, g)$ it is easy to construct sequences of functions $(f_n, g_n)$ such that
$(f_n, g_n) \rightarrow (f, g)$ in $\mathcal{C}^{\alpha} \times
\mathcal{C}^{\beta}$ and such that the pointwise limit $\lim_n f_n g_n$ exists
and is smooth but nonetheless
\[ \lim_n f_n g_n \neq fg. \]
If we want a ``robust'' way to consider the product $fg$ in the situation
$\alpha + \beta < 0$, then we should take this \tmtextit{ambiguity} into
account from the start and think about the product as describing a manifold of
possibilities and not just a single deterministic operation on the inputs.

Part of the analysis of singular SPDEs can be understood as a classification
of these ambiguities: we track the extent to which the possible outcomes of
undefined operations can propagate into the solution theory of a given
equation. We will come back to this point below when we discuss
renormalisations.

Paraproducts and related paradifferential operators have been introduced in
the seminal work of Bony on the propagation of singularities for non-linear
hyperbolic equations~{\cite{bony_calcul_1981,meyer_remarques_1981}}. They
provide a good approximation of non-linear operations, in this case the
product, but can be used also to linearise other operations. For example, the
following paralinearisation result is useful in order to deal with equation
with non-polynomial coefficients: Given a smooth function $f : \mathbb{R}
\rightarrow \mathbb{R}$ we have
\begin{equation}
  u \in \mathcal{C}^{\alpha} \mapsto R_f (u) \assign f (u) - f' (u) \para u
  \in \mathcal{C}^{2 \alpha}, \qquad \alpha > 0, \label{eq:paralinearization}
\end{equation}
which shows that the composition $f (u) \in \mathcal{C}^{\alpha}$ behaves like
the paraproduct $(f' (u) \para u) \in \mathcal{C}^{\alpha}$ modulo a smoother
correction term.

One main tool in the paracontrolled analysis is the following commutator
lemma which describes the interaction between the resonant product and the
paraproduct. For a proof see~{\cite{gubinelli_paracontrolled_2015}}.

\begin{lemma}
  \label{lemma:commutator}Let $\alpha \in (0, 1), \beta, \gamma \in
  \mathbb{R}$ be exponents such that $\beta + \gamma < 0$ and $\alpha + \beta
  + \gamma > 0$. Then there exists a continuous trilinear map $C :
  \mathcal{C}^{\alpha} \times \mathcal{C}^{\beta} \times \mathcal{C}^{\gamma}
  \rightarrow \mathcal{C}^{\alpha + \beta + \gamma}$ such that if $f, g, h$
  are smooth functions we have
  \[ C (f, g, h) = (f \para g) \reso h - f (g \reso h) . \]
\end{lemma}

\section{Paracontrolled analysis}\label{sec:paracontrolled}

The construction of the solution maps proceeds via perturbative analysis with
respect to a linearised approximation of the equation. We look for an
expansion in terms of regularities of the various terms. In the following we ignore the initial conditions for the equations that we treat and when
considering a solution $w$ to $\mathcal{L}w = v$ we always silently assume
that $w (0) = 0$. Once this is understood the case of general initial
conditions does not add much conceptual (although some technical) difficulty. We also ignore the need for renormalisation at the moment and set the (infinite) constants $C$ appearing in the equations equal to $0$ in the following discussion. In general we omit many technical details, an introduction to paracontrolled distributions that provides more details can be found in the lecture notes~\cite{gubinelli_lectures_2015}.

Among the examples we treat, the simplest is the dynamic $\Phi^4_2$
model~(\ref{eq:SQE}), $\varphi : \mathbb{R}_+ \times \mathbb{T}^2
\rightarrow \mathbb{R}$,
\[ \partial_t \varphi (t, x) = \Delta \varphi (t, x) - \varphi (t, x)^3 + \xi (t, x), \qquad t \geqslant 0,
   x \in \mathbb{T}^2. \]
This equation was first solved by Albeverio and R{\"o}ckner~\cite{albeverio_stochastic_1991} in 1991, a solution that is more related to the tools
we present here was given in 2003 by Da Prato and Debussche~\cite{daprato_strong_2003}, and the link between the two solution concepts was recently understood by R\"ockner, Zhu, and Zhu~\cite{rockner_restricted_2015}. Da Prato and Debussche observed that we can decompose the solution as
a sum of two terms which we postulate to be of increasing regularities:
Setting $\varphi = X + \psi$, the equation takes the form
\[ \mathcal{L}X +\mathcal{L} \psi = - X^3 - 3 X^2 \psi - 3 X \psi^2 - \psi^3 +
   \xi \]
where we use the notation $\mathcal{L}= \partial_t - \Delta$ and we choose $X$
in order to cancel the most irregular term on the right hand side, namely the
additive noise. So if $\mathcal{L}X = \xi$, a simple inhomogeneous linear
equation that can be explicitly solved by convolving $\xi$ against the heat
kernel, we are left with an equation for $\psi$,
\begin{equation}
  \mathcal{L} \psi = - X^3 - 3 X^2 \psi - 3 X \psi^2 - \psi^3 .
  \label{eq:dd-step}
\end{equation}
Using stochastic analysis it is possible to prove that
\[ X \in C_T \mathcal{C}^{0 -} \assign \bigcap_{\delta > 0} C ([0, T],
   \mathcal{C}^{- \delta}), \]
and therefore the right hand side of the equation for $\psi$ is still not well
defined because it features the products $X^2$ and $X^3$ of the distribution
$X$ of negative regularity. For the moment we \tmtextit{assume} that $X^3,
X^2$ exist as elements of $C_T \mathcal{C}^{0 -}$. They will be part of the
non-linear features needed to define the solution map. If we also assume that
\[ \psi \in C_T \mathcal{C}^{0 +} \assign \bigcup_{\delta > 0} C ([0, T],
   \mathcal{C}^{\delta}), \]
then the right hand side of eq.~(\ref{eq:dd-step}) is well defined and the
estimates for the product give us $X^2 \psi, X \psi^2 \in C_T \mathcal{C}^{0
-}$. By standard estimates for the heat semigroup generated by the Laplace
operator $\Delta$ the map that sends $v$ to the solution $w$ of $\mathcal{L}w
= v$ is continuous from $C_T \mathcal{C}^{\alpha} = C ([0, T],
\mathcal{C}^{\alpha})$ equipped with $\| v \|_{C_T \mathcal{C}^{\alpha}} =
\sup_{t \in [0, T]} \| v (t) \|_{\mathcal{C}^{\alpha}}$ to $C_T
\mathcal{C}^{\alpha + 2}$. So it follows that we can pose the
equation~(\ref{eq:dd-step}) for $\psi$ as a standard PDE with unique weak
solution $\psi \in C_T \mathcal{C}^{2 -}$ for some\footnote{Possibly $T$ is quite small because when setting up the Picard iteration we pick up a superlinear estimate so at this point we cannot exclude the possibility that the solution blows up in finite time. To prove existence for all times we have to make use of the sign of the nonlinearity $-\varphi^3$ in (\ref{eq:SQE}), see~\cite{mourrat_global_2015}.} $T > 0$. Here the change of variable has
therefore been quite simple: The enhancement $J$, the solution map $\Phi$ and
the projection $\Pi$ take the form
\[ J : \xi \mapsto (X, X^2, X^3), \]
\[ \Phi : \Xi = (\zeta_1, \zeta_2, \zeta_3) \mapsto U = (\zeta_1, \psi),
   \qquad \Pi : U = (\zeta_1, \psi) \mapsto \varphi = \zeta_1 + \psi, \]
where $\psi$ is a weak solution to
\begin{equation}
  \LL \psi = \zeta_3 + 3 \zeta_2 \psi + 3 \zeta_1 \psi^2 + \psi^3 .
  \label{eq:dd-step-abs}
\end{equation}
Note that $\xi$ can be recovered from $(\zeta_1, \zeta_2, \zeta_3) = J (\xi)$
via $\xi =\mathcal{L}X$. Observe also that if $\xi$ is a smooth function, the
composition $\Pi \circ \Phi \circ J : \xi \mapsto \varphi$ gives back the
classical weak solution to eq.~(\ref{eq:SQE}). We will come back again later
to this equation to discuss the constuction of the enhanced data $\Xi$ which
will require a renormalisation in the case of white noise. For details on the
analysis that we sketched above see \cite{daprato_strong_2003}.

While the method we just discussed is simple and elegant, the other singular
equations that we mentioned in the introduction apart from the dynamic
$\Phi^4_2$ model cannot be handled by a simple additive change of variables.
Consider for example the generalised stochastic Burgers
equation~(\ref{eq:gSBE}). Since the noise is additive as in the dynamic
$\Phi^4_2$ equation, we can proceed with the same decomposition of the
solution. We let $X$ be the solution to $\mathcal{L}X = \xi$ and write $u = X
+ v$, where $v$ solves
\begin{equation}
  \mathcal{L}v = G (X + v) \partial_x X + G (X + v) \partial_x v.
  \label{eq:gsbe-transformed}
\end{equation}
The analysis of the regularity of $X$ now gives $X \in C_T \mathcal{C}^{1 / 2
-}$ (it is better behaved than before because we are in one space-dimension
and the white noise becomes more and more irregular in higher dimensions) and
therefore the best we can hope for the right hand side
of~(\ref{eq:gsbe-transformed}) is that it takes values in $C_T \mathcal{C}^{-
1 / 2 -}$ (the regularity of $\partial_x X$) which would put $v$ in $C_T
\mathcal{C}^{3 / 2 -}$, that is $v$ has two degrees of regularity more than
the right hand side which follows from the regularising effect of the
inversion of $\mathcal{L}$ that we discussed above. This would in turn mean
that $G (X + v) \in C_T \mathcal{C}^{1 / 2 -}$ provided $G$ is at least $C^1$.
In this setting the product $G (X + v) \partial_x v$ is well defined since
$\partial_x v$ has regularity $C_T \mathcal{C}^{1 / 2 -}$, but $G (X + v)
\partial_x X$ is not since $\partial_x X$ has negative regularity $C_T
\mathcal{C}^{- 1 / 2 -}$ and we barely fail to ensure that the sum of the
regularities is positive. But no other simple additive subtraction is
available and therefore we need to understand better the structure of the
problematic product in order to determine sufficient conditions to control it.
The paraproduct decomposition gives \
\[ G (X + v) \partial_x X = \underbrace{G (X + v) \para \partial_x X}_{C_T
   \mathcal{C}^{- 1 / 2 -}} + \underbrace{G (X + v) \reso \partial_x X}_{!!} +
   \underbrace{G (X + v) \lpara \partial_x X}_{C_T \mathcal{C}^{0 -}}, \]
where the underbraces denote the respective regularities of the two
paraproducts and the difficulty is isolated in the resonant term. The
paralinearisation result~(\ref{eq:paralinearization}) applied to $G (X + v)$
shows that
\[ G (X + v) = \underbrace{G' (X + v) \para (X + v)}_{C_T \mathcal{C}^{1 / 2
   -}} + \underbrace{R_G (X + v)}_{C_T \mathcal{C}^{1 -}}, \]
and we can decompose the paraproduct on the right hand side into
\[ G' (X + v) \para (X + v) = \underbrace{G' (X + v) \para X}_{C_T
   \mathcal{C}^{1 / 2 -}} + \underbrace{G' (X + v) \para v}_{C_T
   \mathcal{C}^{3 / 2 -}}, \]
which shows that the irregularity of $G (X + v)$ comes from the paraproduct
${G' (X + v) \para X}$ and we can further isolate the difficulty in the resonant
product:
\[ G (X + v) \reso \partial_x X = \underbrace{(G' (X + v) \para X) \reso
   \partial_x X}_{!!} + \underbrace{G' (X + v) \para v}_{C_T \mathcal{C}^{1
   -}} + \underbrace{R_G (X + v) \reso \partial_x X}_{C_T \mathcal{C}^{1 / 2
   -}}, \]
where the last two terms on the right hand side can be controlled by the
estimates for the resonant product because here the sum of the regularities is
strictly positive. To deal with the remaining ill-defined resonant product we
apply the commutator estimate from Lemma~\ref{lemma:commutator} which gives us
\[ (G' (X + v) \para X) \reso \partial_x X = \underbrace{C (G' (X + v), X,
   \partial_x X)}_{C_T \mathcal{C}^{1 / 2 -}} + G' (X + v) \underbrace{(X
   \reso \partial_x X)}_{!!}, \]
where we continue to denote with the underbrace ``!!'' the term which is still
problematic according to the deterministic a priori regularities. The reader
can convince herself that all the other terms are indeed well defined. \ This
analysis allowed us to \tmtextit{isolate} the singular nature of the product
into some non-linear feature pertaining only to the data $X$. Much like the
simpler algebraic analysis allowed by the $\Phi^4_2$ model. Similarly we will
now \tmtextit{assume} a prescribed natural regularity for $X \reso \partial_x
X$, namely $X \reso \partial_x X \in C_T \mathcal{C}^{0 -}$, and include this
term in the enhanced data for this problem. Then the remaining product $G' (X
+ v) (X \reso \partial_x X)$ is well defined because $G' (X + v) \in C_T
\mathcal{C}^{1 / 2 -}$ and therefore the sum of the regularities of the
factors is strictly positive. From here we can continue to solve the equation
for $v$ by a Picard iteration. The enhancement $J$, the solution map $\Phi$
and the projection $\Pi$ now take the form
\[ J : \xi \mapsto (X, X \reso \partial_x X) \]
\[ \Phi : \Xi = (\zeta_1, \zeta_2) \mapsto U = (\zeta_1, v), \qquad \Pi : U =
   (\zeta_1, v) \mapsto u = \zeta_1 + v, \]
where $v \in C_T \mathcal{C}^{3 / 2 -}$ is the solution to the PDE
\[ \LL v = G (\zeta_1 + v) \para \partial_x \zeta_1 + G' (\zeta_1 + v) \zeta_2
   + F (\zeta_1, v) \]
and where $F (\zeta_1, v)$ is a suitable continuous function taking values in
$C_T \mathcal{C}^{0 -}$. For details on this equation we refer to~\cite{hairer_rough_2011} and~\cite{gubinelli_paracontrolled_2015}.

A further level of complexity is exemplified by the gPAM~(\ref{eq:gPAM}). In
this case it is not even possible to start the analysis by an additive change
of variables. The two dimensional space white noise has regularity
$\mathcal{C}^{- 1 -}$, so the best regularity we can hope for $v$ is $v \in
C_T \mathcal{C}^{1 -}$ and then the non-linear term $G (u) \xi$ is not well
defined. The paraproduct decomposition gives
\[ G (u) \xi = \underbrace{G (u) \para \xi}_{C_T \mathcal{C}^{- 1 -}} +
   \underbrace{G (u) \reso \xi}_{!!} + \underbrace{G (u) \lpara \xi}_{C_T
   \mathcal{C}^{0 -}} \]
and proceeding by paralinearisation and commutation we obtain the following
decomposition of the resonant term
\begin{gather}
  G (u) \reso \xi = (G' (u) \para u) \reso \xi + R_G (u) \reso \xi\\
  = \underbrace{C (G' (u), u, \xi)}_{C_T \mathcal{C}^{1 -}} + G' (u)
  \underbrace{(u \reso \xi)}_{!!} + \underbrace{R_G (u) \reso \xi}_{C_T
  \mathcal{C}^{1 -}},
\end{gather}
where we note that $u \reso \xi$ is still not well defined but if we assume it
has its natural regularity $u \reso \xi \in C_T \mathcal{C}^{0 -}$, then the
product $G' (u) (u \reso \xi)$ poses no problem. This means that we can
control the product $G (u) \xi$ once we have a control of the resonant term $u
\reso \xi$. Contrary to the simpler analysis of the gSBE this term is still
quite complex since involves the unknown $u$ and cannot be ``postulated'' as
we did with $X \reso \partial_x X$ before. However, our analysis shows that
the right hand side of the equation can be decomposed in a series of terms of
different regularities, where the worst is $G (u) \para \xi$ (assuming for $u
\reso \xi$ a better regularity). So the solution should satisfy
\[ \mathcal{L}u = G (u) \para \xi + \cdots, \]
where we neglected more regular terms. The idea is then to make a change of
variables to remove this irregular term in the right hand side. A natural
approach is to look for $u$ with a similar form as the right hand side of the
equation, namely a paraproduct plus a smoother remainder, $u = v \para X +
v^{\sharp}$, where $v \in C_T \mathcal{C}^{1 -}, X \in C_T \mathcal{C}^{1 -},
v^{\sharp} \in C_T \mathcal{C}^{2 -}$ are functions to be dermined with
$v^{\sharp}$ more regular than $X$ and $u$. In order to perform this change of
variables in the equation we need to modify the paraproduct $\ast \para \ast$
a bit by introducing some time-smoothing and defining a modified paraproduct
$\ast \mpara \ast$ in terms of which we make the \tmtextit{paracontrolled
Ansatz}:
\begin{equation}
  u = v \mpara X + v^{\sharp} . \label{eq:ansatz}
\end{equation}
This modification of the paraproduct is a small technical point which is not very relevant to the overall
picture. The essential property of $\ast \mpara \ast$ is that the operator
\[ (f, g) \mapsto H (f, g) =\mathcal{L} (f \mpara g) - f \para \mathcal{L}g
\]
maps\footnote{This is not exactly true, we also need some time regularity of $f$ but for simplicity we omit this in the discussion.} the space $C_T \mathcal{C}^{1 -} \times C_T \mathcal{C}^{1 -}$ to $C_T
\mathcal{C}^{0 -}$ despite the fact
that both summands only live in $C_T \mathcal{C}^{- 1 -}$. For the usual
paraproduct $\ast \para \ast$ this is not possible because if we expand
$\mathcal{L} (f \para g) - f \para \mathcal{L}g$ using Leibniz's rule we pick
up the term $\partial_t f \para g$ which cannot be controlled in terms of the
$C_T \mathcal{C}^{1 -}$ norm of $f$. The modified paraproduct overcomes this
difficulty, but there exist other solutions: either use space-time parabolic Besov spaces and related paraproducts (for which the commutator of paraproduct and $\mathcal{L}$
can be controlled by standard estimates), or
define a paraproduct which intertwines exactly with the heat kernel so $H (f,
g) = 0$ (we will discuss this last strategy in more detail in
Section~\ref{sec:higherorder}). In the following we will mostly neglect the
difference between these two paraproducts and always consider $H (f, g)$ as a
smoother remainder term. With this proviso we can expand both sides
of the equation and get
\[ \mathcal{L}u = v \para \mathcal{L}X + H (v, X) +\mathcal{L}v^{\sharp} = G
   (u) \para \xi + \tilde{F} (u, \xi), \]
where $\tilde{F} (u, \xi)$ denotes terms that \tmtextit{should be} in $C_T
\mathcal{C}^{0 -}$. Choosing $v = G (u)$ and $\mathcal{L}X = \xi$ we can get
rid of the irregular term $G (u) \para \xi$ on the right hand side and obtain
an equation for $v^{\sharp}$ which sets it in the good space $C_T
\mathcal{C}^{2 -}$. Now that we have a better description of the solution $u$
via the paracontrolled Ansatz we can go back to the analysis of the resonant
term $u \reso \xi$ and observe that the commutator lemma gives
\[ \underbrace{u \reso \xi}_{!!} = \underbrace{(G (u) \para X) \reso \xi}_{!!}
   + \underbrace{v^{\sharp} \reso \xi}_{C_T \mathcal{C}^{1 -}} = G (u)
   \underbrace{(X \reso \xi)}_{!!} + \underbrace{C (G (u), X, \xi)}_{C_T
   \mathcal{C}^{1 -}} + \underbrace{v^{\sharp} \reso \xi}_{C_T \mathcal{C}^{1
   -}} \]
which again reduces the well-posedness of the right hand side to that of a
polynomial non-linear feature constructed from $\xi$, in this case $X \reso
\xi$. We will assume that $X \reso \xi$ is given as an element of $C_T
\mathcal{C}^{0 -}$ so in particular the product $G (u) (X \reso \xi)$ is then
well defined and in $C_T \mathcal{C}^{0 -}$. Resuming this analysis we can
conclude that the enhancement $J$, the solution map $\Phi$ and the projection
$\Pi$ take here the form
\begin{equation}
  \begin{array}{l}
    J : \xi \mapsto (X, X \reso \xi)\\
    \Phi : \Xi = (\zeta_1, \zeta_2) \mapsto U = (\zeta_1, u, v^{\sharp}),
    \qquad \Pi : U = (\zeta_1, u, v^{\sharp}) \mapsto u,
  \end{array} \label{eq:gPAM-sol}
\end{equation}
where $(u, v^{\sharp})$ is a solution to the system
\[ \left\{ \begin{array}{l}
     \LL v^{\sharp} = G' (u) G (u) \zeta_2 + F (u, v^{\sharp}, \zeta_1)\\
     u = G (u) \mpara \zeta_1 + v^{\sharp}
   \end{array} \right. \]
where $F (u, v^{\sharp}, \zeta_1)$ takes values in $C_T \mathcal{C}^{0 -}$.
The equation for $ (u, v^{\sharp})$, albeit not a standard PDE, can
nonetheless be solved by usual fixpoint methods, at least locally in time\footnote{
We pick up a quadratic estimate from the paralinearisation result~\eqref{eq:paralinearization} because it is based on a second order Taylor expansion, and therefore we cannot exclude that the solutions blows up in finite time. But given an a priori bound on the $L^\infty$ norm of $u$ one can show that $(u,v^\sharp)$ stays bounded in $C_T\CC^{1-} \times C_T \CC^{2-}$, see~\cite{gubinelli_paracontrolled_2015}, and such an a priori bound can for certain nonlinearities $G$ be derived from the maximum principle, see~\cite{cannizzaro_malliavin_2017}.}.

This last example allowed us to introduce the paracontrolled ansatz and the
use of paraproducts to describe spaces of distributions with specific
behaviour. All the other examples of singular SPDEs which we mentioned in the
introduction can be analysed via a change of variables involving linear
combinations of paraproducts and smooth remainder terms. We illustrate the
final result of the analysis instead of going step by step as we did so far.
For the dynamic $\Phi^4_3$ model~(\ref{eq:SQE}) we proceed as for $\Phi^4_2$
and introduce further terms by writing
\[ \varphi = X + Y + \varphi^Q, \qquad \varphi^Q = \psi \mpara Q +
   \varphi^{\sharp}, \]
where the functions $X \in C_T \mathcal{C}^{- 1 / 2 -}, Y \in C_T
\mathcal{C}^{1 / 2 -}, Q \in C_T \mathcal{C}^{1 -}, \psi \in C_T
\mathcal{C}^{1 / 2 -}, \varphi^{\sharp} \in C_T \mathcal{C}^{3 / 2 -}$ solve
\[
   \mathcal{L}X = \xi, \qquad \mathcal{L}Y = - X^3, \qquad \mathcal{L}Q = -X^2,
\]
\begin{equation}
  \mathcal{L} \varphi^{\sharp} = - 3 X^2 \reso Y - 3 \psi (X^2 \reso Q) + F
  (\psi, \varphi^{\sharp}, X, Y, Q), \qquad \psi = 3 (Y + \varphi^Q),
  \label{eq:split-sqe-3}
\end{equation}
where $F$ is a continuous function mapping into $C_T \mathcal{C}^{- 1 / 2 -}$.
As before the main goal of this decomposition is to rewrite all the
problematic resonant products in terms of simple expressions of the driving
noise $\xi$. The enhancement $J$, the solution map $\Phi$ and the projection
$\Pi$ take here the form
\begin{equation}
  \begin{array}{l}
    J : \xi \mapsto (X, Y, Q, X^2 \reso Y, X^2 \reso Q)\\
    \Phi : \Xi = (\zeta_1, \ldots, \zeta_5) \mapsto U = (\zeta_1, \zeta_2,
    \zeta_3, \psi, \varphi^{\sharp}),\\
    \Pi : U = (\zeta_1, \zeta_2, \zeta_3, \varphi, \psi, \varphi^{\sharp})
    \mapsto \zeta_1 + \zeta_2 + \psi \mpara \zeta_3 + \varphi^{\sharp},
  \end{array} \label{eq:SQE-sol}
\end{equation}
where the pair $(\psi, \varphi^{\sharp})$ solves the
equations~(\ref{eq:split-sqe-3}) above with the driving features $(X, Y, Q,
X^2 \reso Y, X^2 \reso Q)$ replaced by generic functions $\Xi = (\zeta_1,
\ldots, \zeta_5)$ with specific regularities, which in this case can be taken
as
\[ \Xi \in C_T \mathcal{C}^{- 1 / 2 -} \times C_T \mathcal{C}^{1 / 2 -} \times
   C_T \mathcal{C}^{1 -} \times C_T \mathcal{C}^{- 1 / 2 -} \times C_T
   \mathcal{C}^{- 0 -} . \]
   The details can be found in~\cite{catellier_paracontrolled_2013}, see also~\cite{mourrat_global_2016} for a proof that solutions exist for all times.
   
In the case of the KPZ equation~(\ref{eq:KPZ}) the decomposition is even more
involved and the enhancement $J$, the solution map $\Phi$ and the projection
$\Pi$ take the form
\begin{gather*}
   J : \xi \mapsto (Y, Y^\zzone, Y^\zztwo, Y^\zzthreereso, Y^\zzfour, \partial_x Y \reso \partial_x
   \mathcal{J}\! Y) \\
   \Phi : \Xi = (\zeta_1, \ldots, \zeta_6) \mapsto U = (\zeta_1, \zeta_2,
   \zeta_3, \psi, h^{\sharp}), \\
   \Pi : U = (\zeta_1, \zeta_2, \zeta_3, \psi, h^{\sharp})
   \mapsto h = \zeta_1 + \zeta_2 + 2 \zeta_3 + \psi \mpara \mathcal{J}\!
   \zeta_1 + h^{\sharp},
\end{gather*}
where $\mathcal{J} v(t) = \int_0^t P_{t-s} v(s) \dd s$ and we recall that $(P_t)_{t\ge0}$ is the heat semigroup, and where
\begin{gather*}
   \LL Y = \xi, \quad \LL Y^\zzone = (\partial_x Y)^2, \quad \LL Y^\zztwo =
   \partial_x Y^\zzone \partial_x Y, \quad \LL Y^\zzthreereso = \partial_x Y^\zztwo \reso \partial_x Y,  \\
   \LL Y^\zzfour = (\partial_x Y^\zzone)^2, \quad \psi = 2 (\psi \mpara \mathcal{J}\! \zeta_1 + h^{\sharp}) + 4 Y^\zztwo, \quad \LL h^{\sharp} = F (\Xi, \psi, h^{\sharp})
\end{gather*}
for a continuous function $F$, see~{\cite{gubinelli_kpz_2017}}. We provide more details for the
CSBE equation~(\ref{eq:CSBE})
\begin{equation}
  \LL u = \chi \partial_x u^2 + \partial_x \xi, \label{eq:CSBE-chi}
\end{equation}
with a general constant $\chi$ in front of the nonlinearity because this will be
needed in Section~\ref{sec:universality}. The change of variables reads
\begin{equation}\label{eq:CSBE-chi-exp}
  u = X + \chi X^{\zzone} + 2 \chi^2 X^{\zztwo} + u^Q, \qquad u' = 2 u^Q + 4
  \chi^2 X^{\zztwo}, \quad u^Q = u' \mpara Q + u^{\sharp},
\end{equation}
where $u^{\sharp}$ solves
\begin{align*}
  \LL u^{\sharp} & = \chi \partial_x u^2 - \chi \partial_x X^2 - 2 \chi^2
  \partial_x (X^{\zzone} X) - \LL (u' \mpara Q)\\
  & = \chi \partial_x (\chi X^{\zzone} + 2 \chi^2 X^{\zztwo} + u^Q)^2 + 2
  \chi \partial_x [X (2 \chi^2 X^{\zztwo} + u^Q)] - \LL (u' \mpara Q)\\
  & = \chi^3 \LL X^{\zzfour} + 2 \chi \partial_x  [u^Q X - u^Q \para X] + 2
  \chi [\partial_x (u^Q  \para X) - u^Q  \para \partial_x X]\\
  & \quad + 4 \chi^3 \LL X^{\zzthreereso} + 4 \chi^3 \partial_x [X^{\zztwo}
  \lpara X] + 4 \chi^3  [\partial_x (X^{\zztwo}  \para X) - X^{\zztwo} \para
  \partial_x X] \\
  &\quad  + \chi \partial_x  [2 \chi X^{\zzone} (2 \chi^2 X^{\zztwo} + u^Q) + (2 \chi^2 X^{\zztwo} + u^Q)^2]  - [ \LL (u' \mpara Q) - u' \para ( \LL Q ) ],
\end{align*}
and the enhanced noise is defined by
\[ \Xi = J (\xi) = (X, X, X^{\zzone}, X^{\zztwo}, X^{\zzthreereso},
   X^{\zzfour}, Q, Q \reso X), \]
where $\LL X = \partial_x \xi,$ $\LL Q = \partial_x X$, and
\[ \quad \quad \LL X^{\zzone} = \partial_x X^2, \quad \LL X^{\zztwo} =
   \partial_x (X X^{\zzone}), \quad \LL X^{\zzthreereso} = \partial_x (X^\zztwo \reso
   X), \quad \LL X^{\zzfour} = \partial_x (X^{\zzone})^2 . \]

\section{Ambiguities and renormalisation}\label{sec:renormalization}

The previous analysis reduces the study of a singular SPDE to that of the
enhancement $J$ and of the enhanced solution map $\Phi$. The enhanced version
of the singular equation is a standard PDE for a new unknown together with
some paradifferential relations. This factorisation ``distillates'' in the
definition of the enhancement $\Xi = J (\xi)$ all the problematic products (or
resonant products). One cannot expect to be able to analyse in full generality
the map $J$ without using specific properties of the driving function $\xi$.
Two basic difficulties are:
\begin{enumerate}
  \item the lack of continuity of the resonant products implies that these
  products are essentially undefined and the enhancement map $J$ can be
  extended to irregular inputs $\xi$ only within a specific approximation
  procedure (or not at all);
  
  \item for most of the above examples, even within the more restricted
  context where we only try to extend $J$ to a given stochastic process $\xi$ through a specific approximation procedure, the enhancement map $J$ fails to extend due to divergences: the
  only natural limit of the resonant products is infinite.
\end{enumerate}
The first difficulty means that any resonant product which is formed without
sufficient regularity should be considered inherently \tmtextit{ambiguous},
that is non-robust to different approximation procedures. A satisfying
analysis of such ambiguities is still lacking in the paracontrolled setting
and much more developed in the setting of regularity structures~{\cite{bruned_algebraic_2016}} where deep connection with algebra and renormalization procedures in Quantum Field Theory have been pointed out.
For some details we suggest the reader to have a look at L.~Zambotti's
contribution in this volume. As an example we take here the case of the
solution theory for the gPAM~(\ref{eq:gPAM-sol}). In order to describe the
effects of the ambiguity on the equation we construct an extended solution map
$\Psi_{\text{ext}}$ for smooth inputs $\eta$ by translations $T_{C} J
(\eta) = (\eta, (\mathcal{J}\! \eta) \reso \eta + C)$ of the enhanced gPAM
noise $J (\eta)$ (here we denote as usual with $\mathcal{J}\! \eta$ a suitable
solution to the parabolic problem $\LL \!\!\! \mathcal{J}\! \eta = \eta$). \ Setting
\[
   \hat{u} = \Psi_{\text{ext}} (\xi, C) \assign \Pi \circ \Phi (T_{C} J(\eta)),
\]
the reader can easily check that the function $\hat{u}$ satisfies a
\tmtextit{modified} PDE which reads
\begin{equation}\label{eq:gPAM-mod}
  \LL \hat{u} = G (\hat{u}) \eta + G' (\hat{u}) G (\hat{u}) C .
\end{equation}
But given $\eta \in C^\infty(\mathbb{T}^2)$ and $C>0$ one can find a sequence of smooth functions $(\eta_n) \subset  C^\infty(\mathbb{T}^2)$ such that $\eta_n$ converges to $\eta$ in $\CC^{-1-}$ but $J(\eta_n)$ converges to $T_CJ(\eta)$. So the analysis of the gPAM model~(\ref{eq:gPAM}) together with the requirement
of stability under perturbations which are only small in very weak topologies
\tmtextit{generates} quite naturally a relaxed equation of the
form~(\ref{eq:gPAM-mod}). Similar considerations are applicable, at least in
principle, to all the other singular SPDEs.

The problem of renormalisation is related to this ambiguity and to the
robustness of the form of the equations under irregular perturbations. Here
the problem is however that certain products are intrinsically
impossible to define due to the presence of infinities and that some
subtraction is required for them to have a finite limit. One of the simplest
situations is still that of the 2d gPAM driven by space white noise. In this
case the product $X \reso \xi$, understood as limit of smooth convolutional
approximations is almost surely infinite. So if want the equations with mollified noise to have a well
defined limit, then we need to start with a suitably modified equation of the form
given by eq.~(\ref{eq:gPAM-mod}) and conceived in such a way that the
additional term provides the necessary cancellations to remove the divergence
in the resonant product. Denoting by $\xi_{\varepsilon}$ the convolutional
regularisation of the white noise $\xi$ and by $c_{\varepsilon}$ a family of
constants we see that if $\hat{u}_{\varepsilon}$ is the solution to
\[
   \LL \hat{u}_{\varepsilon} = G (\hat{u}_{\varepsilon}) \xi_{\varepsilon} - G' (\hat{u}_{\varepsilon}) G (\hat{u}_{\varepsilon}) c_{\varepsilon},
\]
then $\hat{u}_{\varepsilon} = \Psi_{\text{ext}} (\xi_{\varepsilon}, -
c_{\varepsilon}) = \Pi \circ \Phi (T_{- c_{\varepsilon}} J (\eta)) = \Pi \circ
\Phi ((\xi_{\varepsilon}, \mathcal{J}\! \xi_{\varepsilon} \reso
\xi_{\varepsilon} - c_{\varepsilon}))$. In other words $\hat{u}_{\varepsilon}$
is a continuous function of the quantities $\Xi_{\varepsilon} \assign
(\xi_{\varepsilon}, \mathcal{J}\! \xi_{\varepsilon} \reso \xi_{\varepsilon} -
c_{\varepsilon})$. A similar result is also true for the right hand side of the
equation, namely
\[ G (\hat{u}_{\varepsilon}) \diamond \Xi_{\varepsilon} \assign G
   (\hat{u}_{\varepsilon}) \xi_{\varepsilon} - G' (\hat{u}_{\varepsilon}) G
   (\hat{u}_{\varepsilon}) c_{\varepsilon} . \]
By a probabilistic analysis it can be checked that a (non-unique) choice of
$(c_{\varepsilon})_{\varepsilon}$ exists for which
$(\Xi_{\varepsilon})_{\varepsilon}$ converges (in probability and in $L^p$
with respect to the randomness) in the appropriate topology. Denoting by $\Xi$
the limit we have that also the solution converges $\hat{u}_{\varepsilon}
\rightarrow \hat{u}$ and satisfies the \tmtextit{renormalised} singular SPDE
\[ \LL \hat{u} = G (\hat{u}) \diamond \Xi, \]
where the right hand side now is a certain non-linear function of $\hat{u}$ and $\Xi$ which
can be identified with the limit $G (\hat{u}) \diamond \Xi = \lim_{\varepsilon
\rightarrow 0} G (\hat{u}_{\varepsilon}) \diamond \Xi_{\varepsilon}$. In this
particular case the limit is controlled via the paracontrolled
Ansatz for $\hat{u}_{\varepsilon}$ and via the continuity results
which follow from it and from the convergence of the enhanced noise
$\Xi_{\varepsilon}$.

More complex renormalisations are necessary in other equations. For example,
in the case of the $\Phi^4_3$ model the ill-defined products are
\[ X^2, X^3, X^2 \reso \mathcal{J}\! (X^3), X^2 \reso \mathcal{J}\! (X^2), \]
see eq.~(\ref{eq:SQE-sol}), where $X =\mathcal{J}\! (\xi)$. Stochastic analysis
shows that there exists a choice of constants $(c_{\varepsilon}^1,
c^2_{\varepsilon})_{\varepsilon > 0}$ such that $c_{\varepsilon}^1,
c_{\varepsilon}^2 \rightarrow \infty$ as $\varepsilon \rightarrow 0$ for which
the convergence in probability
\[ X^{\diamond 3}_{\varepsilon} \assign X^3_{\varepsilon} - 3
   c_{\varepsilon}^1 X_{\varepsilon} \rightarrow X^{\diamond 3}, \qquad
   X^{\diamond 2}_{\varepsilon} \assign X^2_{\varepsilon} - c_{\varepsilon}^1
   \rightarrow X^{\diamond 2}, \]
\[ X^{\diamond 2}_{\varepsilon} \diamond_{\varepsilon} \mathcal{J}\!
   (X^{\diamond 3}_{\varepsilon}) \assign X^{\diamond 2}_{\varepsilon} \reso
   \mathcal{J}\! (X^{\diamond 3}_{\varepsilon}) - 3 c_{\varepsilon}^2
   X_{\varepsilon} \rightarrow X^{\diamond 2} \diamond \mathcal{J}\!
   (X^{\diamond 3}), \qquad \]
\[ X^{\diamond 2}_{\varepsilon} \diamond_{\varepsilon} \mathcal{J}\!
   (X^{\diamond 2}_{\varepsilon}) \assign X^{\diamond 2}_{\varepsilon} \reso
   \mathcal{J}\! (X^{\diamond 2}_{\varepsilon}) - c_{\varepsilon}^2 \rightarrow
   X^{\diamond 2} \diamond \mathcal{J}\! (X^{\diamond 2}) \]
holds in the appropriate spaces, where $X_{\varepsilon} = \rho_{\varepsilon}
\ast X$ is a convolutional regularisation of $X$ with a smooth kernel. As a
consequence of these convergence results and of the structure of the solution
theory described by eq.~(\ref{eq:SQE-sol}) the function $\varphi_{\varepsilon}
= \Pi \reso \Phi (\Xi_{\varepsilon})$, where
\[ \Xi_{\varepsilon} = (X_{\varepsilon}, \mathcal{J}\! (X^{\diamond
   3}_{\varepsilon}), \mathcal{J}\! (X^{\diamond 2}_{\varepsilon}), X^{\diamond
   2}_{\varepsilon} \diamond_{\varepsilon} \mathcal{J}\! (X^{\diamond
   3}_{\varepsilon}), X^{\diamond 2}_{\varepsilon} \diamond_{\varepsilon}
   \mathcal{J}\! (X^{\diamond 2}_{\varepsilon})), \]
solves the renormalised $\Phi^4_3$ equation
\[ \LL \varphi_{\varepsilon} = - \varphi_{\varepsilon}^3 + 3(c^1_{\varepsilon} +
   c^2_{\varepsilon}) \varphi_{\varepsilon} + \xi_{\varepsilon} \]
and converges to $\varphi = \Pi \reso \Phi (\Xi)$
where $\Xi \assign \lim_{\varepsilon} \Xi_{\varepsilon}$ which again can be
described in terms of a limiting renormalised singular SPDE. In this review we
will not discuss details of the convergence results $\Xi_{\varepsilon}
\rightarrow \Xi$. These are mostly handled with standard probabilistic
techniques. A very comprehensive convergence theory for the stochastic terms
exists in terms of regularity structures \cite{bruned_algebraic_2016, chandra_analytic_2016} but most of the ideas can
be employed also within the paracontrolled approach, see also \cite{mourrat_construction_2016}.

\section{Higher order expansions}\label{sec:higherorder}

The theory of paracontrolled distributions that we discussed so far is
essentially a first order calculus. For example, in the parabolic Anderson
model (gPAM with $G (u) = u$) we expand the solution in terms of a paraproduct
$u \para X$ plus a smoother remainder $u^{\sharp}$, where $\mathcal{L}X = \xi$
and $X \in C_T \mathcal{C}^{\alpha}$, $u^{\sharp} \in C_T \mathcal{C}^{\alpha}$ for $\alpha = 1-$. Then terms like the resonant product $u^{\sharp} \reso
\xi$ pose no problem because $\xi \in \mathcal{C}^{\alpha - 2}$ and
$u^{\sharp} \in C_T \mathcal{C}^{2 \alpha}$ and the sum of the regularities is $3 \alpha
- 2 > 0$. But what if $\xi$ has worse regularity, so if $\alpha \leqslant 2 / 3$ and therefore $3\alpha - 2 \leqslant 0$? One relevant example is the parabolic Anderson
model in $d = 3$, where the space white noise has regularity strictly less
than $- 3 / 2$ and therefore $\alpha < 1 / 2$. The idea, inspired by (controlled) rough paths, is then to find a higher order expansion of $u$ which
should be of the type
\[ u = \sum_{k = 1}^{n - 1} \sum_{\tau \in \mathcal{I}_k} u^{\tau} \para
   X^{\tau} + u^{\sharp} \]
for some $u^{\sharp} \in C_T\mathcal{C}^{n \alpha}$ and suitable
$(X^{\tau})_{\tau \in \bigcup_k \mathcal{I}_k}$, with $X^{\tau} \in C_T
\mathcal{C}^{k \alpha}$ for all $\tau \in \mathcal{I}_k$. If $(n + 1) \alpha -
2 > 0$, then at least the product $u^{\sharp} \reso \xi$ is well defined and
this gives us some hope to construct a continuous map $(u, (u^\tau)_{\tau}, u^{\sharp}, \xi, (X^\tau)_{\tau}, \dots) \mapsto u \xi$.
But what should the $(X^\tau)$ be and how do we define the map? This is not very
obvious and for quite some time the extension of paracontrolled distributions
to more irregular driving noises remained open. But recently Bailleul and
Bernicot~\cite{bailleul_higher_2016} have made significant progress in that
direction.

To understand their results let us have a look at the parabolic Anderson model
\[ \mathcal{L}u = u \xi \]
for $\xi \in \mathcal{C}^{\alpha - 2}$ with $\alpha \in (1 / 2, 2 / 3)$, so
slightly better than the white noise in $d = 3$ but worse than the white noise
in $d = 2$. We guess the expansion $u = \sum_{k = 1}^2 \sum_{\tau \in
\mathcal{I}_k} u^{\tau} \para X^{\tau} + u^{\sharp}$ with $u^{\sharp} \in C_T
\mathcal{C}^{3 \alpha}$ and $(u^{\tau})$ and $(X^{\tau})$ to be determined.
Then
\[ u \xi = u \para \xi_{} + u \lpara \xi + \underbrace{\left( \sum_{k = 1}^2
   \sum_{\tau \in \mathcal{I}_k} u^{\tau} \para X^{\tau} \right) \reso
   \xi}_{!!} + u^{\sharp} \reso \xi,
\]
where as before we single out the problematic resonant product with the underbrace ``!!''. If we assume that $u^{\tau} \in C_T \mathcal{C}^{\alpha}$ for all $\tau \in \mathcal{I}_1 \cup \mathcal{I}_2$ and
$X^{\tau} \in C_T \mathcal{C}^{k \alpha}$ for all $\tau \in \mathcal{I}_k$,
then for $\tau \in \mathcal{I}_2$ the resonant product $(u^{\tau} \para
X^{\tau}) \reso \xi$ can be controlled with the commutator estimate:
\[ (u^{\tau} \para X^{\tau}) \reso \xi = (\underbrace{C (u^{\tau}, X^{\tau},
   \xi)}_{C_T \mathcal{C}^{4 \alpha - 2}} + u^{\tau} (X^{\tau} \reso \xi)), \]
where we used that $\alpha < 1$ and $4 \alpha - 2 > 0$ to see that the
commutator is bounded, and where we need to assume that $X^{\tau} \reso \xi
\in C_T \mathcal{C}^{3 \alpha - 2}$ is extrinsically given. However, for $\tau
\in \mathcal{I}_1$ the term $C (u^{\tau}, X^{\tau}, \xi)$ is still not
well defined because the sum of the regularities is $3 \alpha - 2 < 0$. To
deal with the resonant product $(u^{\tau} \para X^{\tau}) \reso \xi$ we
therefore assume that $u^{\tau}$ is itself paracontrolled of order $2 \alpha$
for all $\tau \in \mathcal{I}_1$:
\[ u^{\tau} = \sum_{\tau' \in \mathcal{I}_1} u^{\tau, \tau'} \para X^{\tau'} +
   u^{\tau, \sharp} \]
with $u^{\tau, \tau'} \in C_T \mathcal{C}^{\alpha}$ and $u^{\tau, \sharp} \in
C_T \mathcal{C}^{2 \alpha}$. Then
\[ (u^{\tau, \sharp} \para X^{\tau}) \reso \xi = C (u^{\tau, \sharp},
   X^{\tau}, \xi) + u^{\tau, \sharp} (X^{\tau} \reso \xi) . \]
At this point we would like to gain $2 \alpha$ degrees of regularity from
$u^{\tau, \sharp}$ in the commutator to see that it is in $C_T \mathcal{C}^{4
\alpha - 2}$, but this is not possible because $2 \alpha > 1$ and the
commutator estimate Lemma~\ref{lemma:commutator} allows us only to gain less
than one derivative. However, the sum of the regularities of $X^{\tau}$ and
$\xi$ is $2 \alpha - 2 > - 1$, and therefore we can use that $u^{\tau,
\sharp} \in C_T \mathcal{C}^{1 -}$ to obtain
\begin{align*}
   \| C (u^{\tau, \sharp}, X^{\tau}, \xi) \|_{C_T \mathcal{C}^{(1 -) + 2 \alpha - 2}} & \lesssim \| u^{\tau, \sharp} \|_{C_T \mathcal{C}^{1 -}} \| X^{\tau} \|_{C_T \mathcal{C}^{\alpha}} \| \xi \|_{C_T \mathcal{C}^{\alpha -  2}} \\
   & \leqslant \| u^{\tau, \sharp} \|_{C_T \mathcal{C}^{2 \alpha}} \| X^{\tau} \|_{C_T \mathcal{C}^{\alpha}} \| \xi \|_{C_T \mathcal{C}^{\alpha -  2}}.
\end{align*}
Since we estimate the commutator in a space of positive regularity and $3\alpha-2 < 0$, we still get $C (u^{\tau, \sharp}, X^{\tau}, \xi)
\in C_T \mathcal{C}^{3 \alpha - 2}$.

\begin{remark}
  \label{rmk:too-little-gain}This argument is particular to the
  not-so-singular problem studied here and breaks down if $\alpha < 1 / 2$. In
  that case it may be necessary to develop a version of the commutator estimate
  which allows to gain more than one derivative from $u^{\tau, \sharp}$. This
  can be achieved by subtracting not only $u^{\tau, \sharp} (X^{\tau} \reso
  \xi)$ from $(u^{\tau, \sharp} \para X^{\tau}) \reso \xi$ but also further
  correction terms that involve modified Littlewood-Paley blocks
  and roughly speaking correspond to polynomial terms in regularity
  structures. Currently there is no reference where this is worked out.
\end{remark}

Next, we note that the product $u^{\tau, \sharp} (X^{\tau} \reso \xi)$ is
well defined if $X^{\tau} \reso \xi \in C_T \mathcal{C}^{2 \alpha - 2}$ is
given because then the sum of the regularities of its factors is $4 \alpha - 2
> 0$. It remains to understand the resonant product
\[ ((u^{\tau, \tau'} \para X^{\tau'}) \para X^{\tau}) \reso \xi = C (u^{\tau,
   \tau'} \para X^{\tau'}, X^{\tau}, \xi) + (u^{\tau, \tau'} \para X^{\tau'})
   (X^{\tau} \reso \xi) \]
where our commutator estimate really fails: After all $u^{\tau, \tau'} \para X^{\tau'} \in C_T \mathcal{C}^{\alpha}$,
$X^{\tau} \in C_T \mathcal{C}^{\alpha}$, and $\xi \in C_T \mathcal{C}^{\alpha
- 2}$ so the sum of the regularities is $3 \alpha - 2 < 0$. But Bailleul and
Bernicot realised that one can iterate the commutator estimate in the
following way:

\begin{lemma}[\cite{bailleul_higher_2016}, formula(3.8)]  
  Let $\alpha, \beta, \gamma, \delta \in \mathbb{R}$ be exponents such that $\alpha + \beta + \gamma + \delta > 0$, then
  there exists a four-linear map $C^{(2)} : \mathcal{C}^{\alpha} \times
  \mathcal{C}^{\beta} \times \mathcal{C}^{\gamma} \times \mathcal{C}^{\delta}
  \rightarrow \mathcal{C}^{\alpha + \beta + \gamma + \delta}$ such that if $f,
  g, h, \zeta$ are smooth functions we have
  \[ C^{(2)} (f, g, h, \zeta) = C (f \para g, h, \zeta) - f C (g, h, \zeta) .
  \]
\end{lemma}

Therefore, we can set
\[ C (u^{\tau, \tau'} \para X^{\tau'}, X^{\tau}, \xi) = C^{(2)} (u^{\tau,
   \tau'}, X^{\tau'}, X^{\tau}, \xi) + u^{\tau, \tau'} C (X^{\tau'}, X^{\tau},
   \xi) \]
which is well defined provided that $C (X^{\tau'}, X^{\tau}, \xi)$ is given and has its natural
regularity $C_T \mathcal{C}^{3 \alpha - 2}$. The only remaining product that
we still need to control is then
\begin{align*}
  (u^{\tau, \tau'} \para X^{\tau'}) (X^{\tau} \reso \xi) & = (u^{\tau, \tau'}
  \para X^{\tau'}) \para (X^{\tau} \reso \xi) + (u^{\tau, \tau'} \para
  X^{\tau'}) \lpara (X^{\tau} \reso \xi)\\
  & \quad + C (u^{\tau, \tau'}, X^{\tau'}, X^{\tau} \reso \xi) + u^{\tau, \tau'}
  (X^{\tau'} \reso (X^{\tau} \reso \xi)),
\end{align*}
which is under control as long as $X^{\tau'} \reso (X^{\tau} \reso \xi) \in
C_T \mathcal{C}^{3 \alpha - 2}$ is given. So let us put everything together:

\begin{align*}
  u \xi & = \underbrace{u \lpara \xi}_{2 \alpha - 2} + \underbrace{u \para
  \xi}_{\alpha - 2} + \underbrace{u^{\sharp} \reso \xi}_{4 \alpha - 2} +
  \sum_{\tau \in \mathcal{I}_2} (\underbrace{C (u^{\tau}, X^{\tau}, \xi)}_{4
  \alpha - 2} + \underbrace{u^{\tau} (X^{\tau} \reso \xi)}_{3 \alpha - 2})\\
  &\quad + \sum_{\tau \in \mathcal{I}_1} (\underbrace{C (u^{\tau, \sharp}, X^{\tau},
  \xi)}_{3 \alpha - 2} + \underbrace{u^{\tau, \sharp} (X^{\tau} \reso \xi)}_{2
  \alpha - 2})\\
  &\quad + \sum_{\tau, \tau' \in \mathcal{I}_1} (\underbrace{C^{(2)} (u^{\tau,
  \tau'}, X^{\tau'}, X^{\tau}, \xi)}_{4 \alpha - 2} + \underbrace{u^{\tau,
  \tau'} C (X^{\tau'}, X^{\tau}, \xi)}_{3 \alpha - 2})\\
  &\quad + \sum_{\tau, \tau' \in \mathcal{I}_1} (\underbrace{C (u^{\tau, \tau'},
  X^{\tau'}, X^{\tau}, \xi)}_{4 \alpha - 2} + \underbrace{u^{\tau, \tau'}
  (X^{\tau'} \reso (X^{\tau} \reso \xi))}_{3 \alpha - 2})\\
  &\quad + \sum_{\tau, \tau' \in \mathcal{I}_1} (\underbrace{(u^{\tau, \tau'} \para
  X^{\tau'}) \para (X^{\tau} \reso \xi)}_{2 \alpha - 2} +
  \underbrace{(u^{\tau, \tau'} \para X^{\tau'}) \lpara (X^{\tau} \reso \xi)}_{3
  \alpha - 2}),
\end{align*}
where the underbrace indicates the regularity of each term. Therefore, the product $u \xi$ is under control if all of the following terms are extrinsically given
\begin{equation}\label{eq:extended-distribution} \{ X^{\sigma} \reso \xi, X^{\tau} \reso
  \xi, C (X^{\tau'}, X^{\tau}, \xi), X^{\tau'} \reso (X^{\tau} \reso \xi) :
  \tau, \tau' \in \mathcal{I}_1, \sigma \in \mathcal{I}_2 \}
\end{equation}
and have their natural regularity.

But making sense of the product $u \xi$ is only the first step in the analysis
of the equation. Next, we should check that for a paracontrolled $u$ also the
solution $v$ to $\mathcal{L}v = u \xi$ is paracontrolled. Here a problem
occurs: above we discussed that if we make the paracontrolled Ansatz
\[ v = \sum_{k = 1}^2 \sum_{\tau \in \mathcal{I}_k} v^{\tau} \mpara
   X^{\tau} + u^{\sharp} \]
(now with the modified paraproduct $\mpara$ rather than $\para$), then we
can commute the operator $\mathcal{L}$ with $\mpara$ in the sense that
$\mathcal{L} (v^{\tau} \mpara X^{\tau}) - v^{\tau} \para
\mathcal{L}X^{\tau}$ has higher regularity. However, the commutation can gain at most
one degree of regularity from $v^{\tau}$, and as we have just seen in
Remark~\ref{rmk:too-little-gain} this may not always be sufficient. So Bailleul
and Bernicot introduce an ``intertwined'' paraproduct
defined as
\[ (f \olessthan g) (t) = \int_0^t P_{t - s} (f \para \mathcal{L}g) (s) \mathd
   s, \]
where $(P_t)_{t \geqslant 0}$ is the heat semigroup generated by $\Delta$.
Then by definition the exact relation $\mathcal{L} (f \olessthan g) = f \para
\mathcal{L}g$ holds, without error term. So we make the modified
paracontrolled Ansatz
\[ v = \sum_{k = 1}^2 \sum_{\tau \in \mathcal{I}_k} v^{\tau} \olessthan
   X^{\tau} + v^{\sharp} \]
with the same regularities as above. Then
\[ \mathcal{L}v = \sum_{k = 1}^2 \sum_{\tau \in \mathcal{I}_k} v^{\tau} \para
   \mathcal{L}X^{\tau} +\mathcal{L}v^{\sharp}, \]
and if we take $\tau_1 \in \mathcal{I}_1$ with $\mathcal{L}X^{\tau_1} = \xi$
and $v^{\tau_1} = u$, then $v^{\tau_1} \para \mathcal{L} \xi^{\tau}$ cancels
the worst regularity contribution $u \para \xi \in C_T \mathcal{C}^{\alpha}$
to $u \xi$. The remaining $X^{\tau}$ and $v^{\tau}$ have to be chosen such
that all contributions of regularity $2 \alpha - 2$ cancel. The most tricky
term to deal with is $u \lpara \xi = \xi \para u$ which we decompose as
\[ \xi \para u = \xi \para \left( \sum_{k = 1}^2 \sum_{\tau \in \mathcal{I}_k}
   u^{\tau} \para X^{\tau} + u^{\sharp} \right) = \sum_{\tau \in
   \mathcal{I}_1} \xi \para (u^{\tau} \para X^{\tau}) + \sum_{\tau \in
   \mathcal{I}_2} \underbrace{\xi \para (u^{\tau} \para X^{\tau})}_{3 \alpha -
   2} + \underbrace{\xi \para u^{\sharp}}_{4 \alpha - 2} . \]
The first sum on the right hand side still has regularity $2 \alpha - 2$, but
Theorem~8 of~\cite{bailleul_higher_2016} gives
\[ T (\xi, u^{\tau}, X^{\tau}) = \xi \para (u^{\tau} \para X^{\tau}) -
   u^{\tau} \para (\xi \para X^{\tau}) \in C_T \mathcal{C}^{3 \alpha - 2} . \]
%\tmcolor{red}{They have a restriction that $X^{\tau}$ must have negative
%regularity, which seems absurd because the Theorem was written to deal with
%situations like this. Ignore this for now.}
Thus, we have
\[ \xi \para u = \sum_{\tau \in \mathcal{I}_1} \underbrace{T (\xi, u^{\tau},
   X^{\tau})}_{3 \alpha - 2} + \sum_{\tau \in \mathcal{I}_2} \underbrace{\xi
   \para (u^{\tau} \para X^{\tau})}_{3 \alpha - 2} + \underbrace{\xi \para
   u^{\sharp}}_{4 \alpha - 2} + \sum_{\tau \in \mathcal{I}_1}
   \underbrace{u^{\tau} \para (\xi \para X^{\tau})}_{2 \alpha - 2}, \]
and we need a contribution from $\LL v$ to cancel the last term on the right hand
side. The other terms in the product $u \xi$ which have regularity worse than
$3 \alpha - 2$ are $\sum_{\tau \in \mathcal{I}_1} u^{\tau, \sharp} (X^{\tau}
\reso \xi)_{}$ and \ $\sum_{\tau, \tau' \in \mathcal{I}_1} (u^{\tau, \tau'}
\para X^{\tau'}) \para (X^{\tau} \reso \xi)$, and for fixed $\tau \in
\mathcal{I}_1$ we have
\[ u^{\tau, \sharp} (X^{\tau} \reso \xi) - u^{\tau, \sharp} \para (X^{\tau}
   \reso \xi) \in C_T \mathcal{C}^{4 \alpha - 2}, \]
so overall the contribution of regularity $2 \alpha - 2$ is
\[ \sum_{\tau \in \mathcal{I}_1} \left( u^{\tau, \sharp} \para (X^{\tau} \reso
   \xi) + \sum_{\tau' \in \mathcal{I}_1} (u^{\tau, \tau'} \para X^{\tau'})
   \para (X^{\tau} \reso \xi) \right) = \sum_{\tau \in \mathcal{I}_1} u^{\tau}
   \para (X^{\tau} \reso \xi) . \]
In conclusion, we obtain the decomposition
\[ u \xi - \underbrace{u \para \xi}_{\alpha - 2} - \underbrace{\sum_{\tau \in
   \mathcal{I}_1} u^{\tau} \para \{ (\xi \para X^{\tau}) + (X^{\tau} \reso
   \xi) \}}_{2 \alpha - 2} \in C_T \mathcal{C}^{3 \alpha - 2} . \]
We can rewrite $X^{\tau} \reso \xi + \xi \para X^{\tau} = X^{\tau} \xi -
X^{\tau} \para \xi$, and therefore we should set $\mathcal{I}_1 = \{ 1
\}$ with $\mathcal{L}X^{1} = \xi$ and $\mathcal{I}_2 = \{ 2 \}$ with
$\mathcal{L}X^{2} = X^{1} \xi - X^{1} \para \xi$. Then we get
with the Ansatz $v^{1} = u$ and $v^{2} = u^{1}$ that
\[ \mathcal{L}v = u \para \xi + u^{1} \para (X^{1} \xi - X^{1}
   \para \xi) +\mathcal{L}v^{\sharp}, \]
so if we set $\mathcal{L}v = u \xi$, then
\[ \mathcal{L}v^{\sharp} = u \xi - u \para \xi - u^{1} \para (X^{1}
   \xi - X^{1} \para \xi) \in C_T \mathcal{C}^{3 \alpha - 2}, \]
which proves that $v^{\sharp} \in C_T \mathcal{C}^{3 \alpha}$ and therefore
the paracontrolled Ansatz was justified. Moreover, since $v^{1} = u$ we
have with $v^{1, 1} = u^{1}$ that $v^{1} - v^{1,
1} \para X^{1} = u - u^{1} \para X^{1} \in C_T
\mathcal{C}^{2 \alpha}$ and therefore also $v^{1}$ is paracontrolled. The terms in~(\ref{eq:extended-distribution}) that we need to construct in order to make sense of all the products are
\[ X^{1} \reso \xi, X^{2} \reso \xi, C (X^{1}, X^{1},
   \xi), X^{1} \reso (X^{1} \reso \xi) . \]
But of course now we were inconsistent: We started with $u = u^{1} \para
X^{1} + u^{2} \para X^{2} + u^{\sharp}$ and ended up with a $v$
that on the first level is paracontrolled in terms of the new intertwined
paraproduct,
\[ v - v^{1} \olessthan X^{1} - v^{2} \olessthan X^{2} \in
   C_T \mathcal{C}^{3 \alpha}, \]
but on the second level is paracontrolled in terms of the usual paraproduct,
$v^{1} - v^{1, 1} \para X^{1} \in C_T \mathcal{C}^{2
\alpha}$. To set up a Picard iteration the map that sends $u$ to the solution
$v$ of $\mathcal{L}v = u \xi$ should map the space of paracontrolled
distributions into itself, so we should assume that also $u$ was
paracontrolled in terms of $\ast\olessthan\ast$ and also $v^{1}$ is
paracontrolled in terms of $\ast\olessthan\ast$. For that we need to understand the
relation between the two paraproducts $\ast \para \ast$ and $\ast \olessthan
\ast$ and also some commutator estimates involving $\ast\olessthan\ast$. All this is
worked out in~\cite{bailleul_higher_2016}, where also the nonlinear case of
gPAM with $\xi$ of regularity $\xi \in \CC^{\alpha-2}$ with $\alpha < 2/3$ is treated. In that case we also need a higher
order version of the paralinearisation result~(\ref{eq:paralinearization}), but this is relatively easy to derive.

\section{Weak universality}\label{sec:universality}

One prominent application of the theory of singular SPDEs is the derivation of
scaling limits of random fields described by local non-linear stochastic
dynamics. As an example we will sketch the case of the CSBE equation which has
been first analysed via regularity structures by Hairer and
Quastel~{\cite{hairer_class_2015}}. As a mesoscopic model of a
weakly-asymmetric diffusion we will consider the solution of the following
SPDE. Take a small $\varepsilon > 0$ and let $\mathbb{T}_{\varepsilon}
=\mathbb{T}/ \varepsilon$ and $v\colon \mathbb{R}_+ \times
\mathbb{T}_{\varepsilon} \rightarrow \mathbb{R}$ be the solution to
\[ \LL v = \varepsilon^{1 / 2} \partial_x P (v) + \partial_x \eta \]
where $\eta$ is a Gaussian noise on $\mathbb{R}_+ \times
\mathbb{T}_{\varepsilon}$ with finite-range space-time correlations and $P\colon
\mathbb{R} \rightarrow \mathbb{R}$ is a given smooth function. We assume that
$\eta$ is centred with covariance
\[ \mathbb{E} [\eta (t, x) \eta (s, y)] = C_{\varepsilon} (t - s, x - y),
   \qquad t, s \in \mathbb{R}, x, y \in \mathbb{T}_{\varepsilon}, \]
where $C_{\varepsilon}$ is (in the second variable) the $\mathbb{T}_{\varepsilon}$-periodised version of
a function $C \colon \mathbb R \times \mathbb R \to \mathbb R$ which has with sufficient polynomial decay in both
variables. The parabolic change of variables $v_{\varepsilon} (t, x) =
\varepsilon^{- 1 / 2} v (t / \varepsilon^2, x / \varepsilon)$ gives the
equation
\[ \LL v_{\varepsilon} = \varepsilon^{- 1} \partial_x P (\varepsilon^{1 / 2}
   v_{\varepsilon}) + \partial_x \xi_{\varepsilon} \]
where $\xi_{\varepsilon} = \varepsilon^{- 3 / 2} \eta (\cdot / \varepsilon^2,
\cdot / \varepsilon)$ is a noise which converges to a space-time white
noise $\xi$. Indeed the rescaled fields $v_{\varepsilon}, \xi_{\varepsilon}$
live on the standard torus $\mathbb{T}$ and
\[ \mathbb{E} [\xi_{\varepsilon} (t, x) \xi_{\varepsilon} (s, y)] =
   \varepsilon^{- 3} C_{\varepsilon} ((t - s) / \varepsilon^2, (x - y) /
   \varepsilon), \qquad t, s \in \mathbb{R}, x, y \in \mathbb{T},
\]
where, as $\varepsilon \rightarrow 0$ we have $\varepsilon^{- 3}
C_{\varepsilon} ((t - s) / \varepsilon^2, (x - y) / \varepsilon) \rightarrow
\delta (t - s) \delta (x - y)$ weakly as a space-time distribution. The goal
of the analysis is to show is that as $\varepsilon \rightarrow 0$ the function
$v_{\varepsilon}$ converges to the solution of the CSBE equation~(\ref{eq:CSBE})
with a specific constant $\chi$ in front of the non-linearity. The constant
will depend only on the shape of the function $P$. By going to a reference
frame via a constant velocity $a_{\varepsilon}$ change of variables we have
that $v_{\varepsilon} = v_{\varepsilon} (t, x) =
\varepsilon^{- 1 / 2} v (t / \varepsilon^2, (x  + a_\varepsilon t) / \varepsilon )$ is the solution to
\[
   \LL v_{\varepsilon} = a_{\varepsilon} \partial_x v_{\varepsilon} +
   \varepsilon^{- 1} \partial_x P (\varepsilon^{1 / 2} v_{\varepsilon}) +
   \partial_x \xi_{\varepsilon} . = \varepsilon^{- 1} \partial_x \tilde{P}
   (\varepsilon^{1 / 2} v_{\varepsilon}) + \partial_x \xi_{\varepsilon},
\]
where $\tilde{P} (x) = P (x) + a_{\varepsilon} \varepsilon^{1 / 2} x$. We
make an Ansatz of the form
\[ v_{\varepsilon} = X_{\varepsilon} + \chi \tilde{X}^{\zzone}_{\varepsilon} +
   2 \chi^2 \tilde{X}^{\zztwo}_{\varepsilon} + v_{\varepsilon}^Q,
\]
where $\chi$ is a real number different from $0$ and $\tilde{X}^{\zzone}_{\varepsilon}$ and $ \tilde{X}^{\zztwo}_{\varepsilon}$ are functions, all to be determined later. We then use a Taylor expansion to get
\begin{align*}
   \frac{1}{\varepsilon} \tilde{P} (\varepsilon^{1 / 2} v_{\varepsilon}) & =
   \frac{1}{\varepsilon} \tilde{P} (\varepsilon^{1 / 2} X_{\varepsilon}) +
   \frac{1}{\varepsilon^{1 / 2}} \tilde{P}' (\varepsilon^{1 / 2}
   X_{\varepsilon}) (\chi \tilde{X}^{\zzone}_{\varepsilon} + 2 \chi^2
   \tilde{X}^{\zztwo}_{\varepsilon} + v_{\varepsilon}^Q) \\
   &\quad + \frac{1}{2} \tilde{P}'' (\varepsilon^{1 / 2} X_{\varepsilon}) (\chi
   \tilde{X}^{\zzone}_{\varepsilon} + 2 \chi^2
   \tilde{X}^{\zztwo}_{\varepsilon} + v_{\varepsilon}^Q)^2 + R_{\varepsilon},
\end{align*}
where
\[ R_{\varepsilon} = \varepsilon^{1 / 2} \frac{1}{2} \int_0^1 \mathd \tau (1-\tau)^2
   \tilde{P}^{(3)} (\varepsilon^{1 / 2} X_{\varepsilon} + \tau \varepsilon^{1
   / 2} (\chi \tilde{X}^{\zzone}_{\varepsilon} + 2 \chi^2
   \tilde{X}^{\zztwo}_{\varepsilon} + v_{\varepsilon}^Q)) (\chi
   \tilde{X}^{\zzone}_{\varepsilon} + 2 \chi^2
   \tilde{X}^{\zztwo}_{\varepsilon} + v_{\varepsilon}^Q)^3 . \]
Let assume that $\tilde{P}^{(3)} (x) = P^{(3)} (x)$ has polynomial growth of
order $M$ and that we have the following bounds
\[ \| X_{\varepsilon} \|_{L^{\infty}} \lesssim \varepsilon^{- 1 / 2 - \kappa},
   \quad \| \chi \tilde{X}^{\zzone}_{\varepsilon} + 2 \chi^2
   \tilde{X}^{\zztwo}_{\varepsilon} + v_{\varepsilon}^Q \|_{L^{\infty}}
   \lesssim \varepsilon^{0 - \kappa} \]
for some small $\kappa > 0$. Then
\[ \| R_{\varepsilon} \|_{L^{\infty}} \lesssim \varepsilon^{1 / 2}
   \varepsilon^{- \kappa M} \varepsilon^{- 3 \kappa} \lesssim \varepsilon^{1 /
   2 - \kappa (M + 3)} \]
so if $\kappa$ is small enough this remainder goes to zero in $L^{\infty}$.
This justifies the Taylor expansion at least under the assumptions we made.

Now we set
\begin{gather*}
   \chi_{\varepsilon} = \tilde{P}'' (\varepsilon^{1 / 2} X_{\varepsilon}), \quad \LL X_{\varepsilon} = \partial_x \xi_{\varepsilon}, \\
   2 \chi \tilde{X}_{\varepsilon} = \frac{1}{\varepsilon^{1 / 2}} \tilde{P}' (\varepsilon^{1 / 2} X_{\varepsilon}), \quad \LL \tilde{Q}_{\varepsilon} = \partial_x \tilde{X}_{\varepsilon}, \quad \chi \LL  \tilde{X}^{\zzone}_{\varepsilon} = \frac{1}{\varepsilon} \partial_x \tilde{P} (\varepsilon^{1 / 2} X_{\varepsilon}),
\end{gather*}
and
\[
   \LL \tilde{X}^{\zztwo}_{\varepsilon} = \partial_x ( \tilde{X}^{\zzone}_{\varepsilon} \reso \tilde{X}_{\varepsilon}), \quad \LL \tilde{X}^{\zzthreereso}_\varepsilon = \partial_x ( \tilde{X}^{\zztwo}_{\varepsilon}\reso \tilde{X}_{\varepsilon} ), \quad \chi \LL \tilde{X}^{\zzfour}_{\varepsilon} = \frac{1}{2} \partial_x [ \chi_{\varepsilon} (\tilde{X}^{\zzone}_{\varepsilon})^2].
 \]
Note that $\chi$ is a real number and $\chi_\varepsilon$ is a function and the two do not agree. The reason for the notation is that in the end $\chi_\varepsilon$ will converge to $\chi$. With these definitions at hand we see that setting
\[ v_{\varepsilon} = X_{\varepsilon} + \chi \tilde{X}^{\zzone}_{\varepsilon} +
   2 \chi^2 \tilde{X}^{\zztwo}_{\varepsilon} + v^Q_{\varepsilon}, \qquad
   v'_{\varepsilon} = 2 \chi v^Q_{\varepsilon} + 4 \chi^3
   \tilde{X}^{\zztwo}_{\varepsilon}, \quad v^Q_{\varepsilon} =
   v'_{\varepsilon} \mpara \tilde{Q}_{\varepsilon} + v^{\sharp}_{\varepsilon} \]
we have
\begingroup
\allowdisplaybreaks
\begin{align*}
  \LL v_{\varepsilon}^{\sharp} & = \frac{1}{\varepsilon} \partial_x [\tilde{P}
  (\varepsilon^{1 / 2} v_{\varepsilon})] - \frac{1}{\varepsilon} \partial_x [\tilde{P} (\varepsilon^{1 / 2} X_{\varepsilon})] - 2 \chi^2 \partial_x [( \tilde{X}^{\zzone}_{\varepsilon}  \tilde{X}_{\varepsilon})] - \LL
  (v_{\varepsilon}' \mpara \tilde{Q}_{\varepsilon})\\
  & = \frac{1}{\varepsilon^{1 / 2}} \partial_x [\tilde{P}' (\varepsilon^{1 /
  2} X_{\varepsilon}) (2 \chi^2 \tilde{X}^{\zztwo}_{\varepsilon} +
  v_{\varepsilon}^Q)] + \frac{1}{2} \partial_x [\tilde{P}'' (\varepsilon^{1 / 2}
  X_{\varepsilon}) (\chi \tilde{X}^{\zzone}_{\varepsilon} + 2 \chi^2
  \tilde{X}^{\zztwo}_{\varepsilon} + v_{\varepsilon}^Q)^2] \\
  &\quad + \partial_x R_{\varepsilon} - \LL (v_{\varepsilon}' \mpara \tilde{Q}_{\varepsilon})\\
  & = \chi 2 \partial_x [\tilde{X}_{\varepsilon} (2 \chi^2
  \tilde{X}^{\zztwo}_{\varepsilon} + v_{\varepsilon}^Q)] + \frac{1}{2}
  \partial_x [\chi_{\varepsilon} (\chi \tilde{X}^{\zzone}_{\varepsilon} + 2
  \chi^2 \tilde{X}^{\zztwo}_{\varepsilon} + v_{\varepsilon}^Q)^2] + \partial_x
  R_{\varepsilon} - \LL (v_{\varepsilon}' \mpara \tilde{Q}_{\varepsilon})\\
  & = \chi^3 \LL \tilde{X}^{\zzfour}_{\varepsilon} + 2 \chi \partial_x 
  [v_{\varepsilon}^Q  \tilde{X}_{\varepsilon} - v_{\varepsilon}^Q \para
  \tilde{X}_{\varepsilon}] + 2 \chi [\partial_x (v_{\varepsilon}^Q \para 
  \tilde{X}_{\varepsilon}) - v_{\varepsilon}^Q \para \partial_x 
  \tilde{X}_{\varepsilon}]\\
  & \quad + 4 \chi^3 \LL \tilde{X}^{\zzthreereso}_{\varepsilon} + 4 \chi^3
  \partial_x [\tilde{X}^{\zztwo}_{\varepsilon} \lpara \tilde{X}_{\varepsilon}]
  + 4 \chi^3  [\partial_x (\tilde{X}^{\zztwo}_{\varepsilon}  \para
  \tilde{X}_{\varepsilon}) - \tilde{X}^{\zztwo}_{\varepsilon}  \para
  \partial_x  \tilde{X}_{\varepsilon}] \\
  &\quad + \partial_x \chi_{\varepsilon} [2 \chi
  \tilde{X}^{\zzone}_{\varepsilon} (2 \chi^2 \tilde{X}^{\zztwo}_{\varepsilon}
  + v_{\varepsilon}^Q) + (2 \chi^2 \tilde{X}^{\zztwo}_{\varepsilon} +
  v_{\varepsilon}^Q)^2]  - \left[ \LL (v_{\varepsilon}' \mpara \tilde{Q}) -
  v_{\varepsilon}' \para ( \LL \tilde{Q} ) \right] \\
  &\quad + \partial_x
  R_{\varepsilon},
\end{align*}
\endgroup
which has to be compared with the expansion that we used for the
solution to the CSBE~(\ref{eq:CSBE-chi-exp}). The only structural difference
is the presence of random fields $\chi_{\varepsilon}$ in the place of some of
the $\chi$ and the additional source term $\partial_x R_{\varepsilon}$ which
however goes to zero in a topology that is compatible with the required regularity of $v^{\sharp}_{\varepsilon}$. The role of the enhancement $\Xi$ is taken by the
(slighly modified) family of random fields
\[ \tilde{\Xi}_{\varepsilon} = (\chi_{\varepsilon}, X_{\varepsilon},
   \tilde{X}_{\varepsilon}, \tilde{X}^{\zzone}_{\varepsilon},
   \tilde{X}^{\zztwo}_{\varepsilon}, \tilde{X}^{\zzthreereso}_{\varepsilon},
   \tilde{X}^{\zzfour}_{\varepsilon}, \tilde{Q}_{\varepsilon}, \tilde{Q}_{\varepsilon} \reso
   \tilde{X}_{\varepsilon}) . \]
If we prove that $\tilde{\Xi}_{\varepsilon}$ converges to an enhancement
$\tilde{\Xi}$
\[ \tilde{\Xi} = (\chi, X, X, X^{\zzone}, X^{\zztwo}, X^{\zzthreereso},
   X^{\zzfour}, Q, Q \reso X), \]
in suitable topologies, where $\chi$ is a constant, we will have automatically that $(v_{\varepsilon},
v'_{\varepsilon}, v^{\sharp}_{\varepsilon}) \rightarrow (u, u', u^{\sharp})$
in the appropriate topologies. It is not our aim here to fully develop this
sketch of proof, especially because a non-trivial part consists
in proving the convergence of the stochastic data for which some powerful machinery
has been devised in the paper by Hairer and Quastel~{\cite{hairer_class_2015}}.

Let us just mention some specific mechanisms at work in the convergence
result. The enhanced data for $v^\varepsilon$ involves non-linear functions of the random field
$\varepsilon^{1 / 2} X_{\varepsilon}$. A direct computation shows that for any
given $(t, x)$ the family of random variables $(\varepsilon^{1 / 2} X_{\varepsilon} (t, x))_\varepsilon$ converges to a Gaussian random variable with a finite variance. In general the covariance of the random field
$X_{\varepsilon}$ is given by
\[ Q_{\varepsilon} (y - z, t - s) =\mathbb{E} [X_{\varepsilon} (t, y)
   X_{\varepsilon} (s, z)] = \varepsilon^{- 1} (e^{\Delta (t - s) /
   \varepsilon^2} C) ((y - z) / \varepsilon). \]
This quantity allows for various bounds:
\[ | Q_{\varepsilon} (y - z, t - s) | \lesssim \varepsilon^{- 1} \wedge (t -
   s)^{- 1 / 2} \wedge (y - z)^{- 1} . \]
Some Gaussian analysis can be used to show that if we have a sequence of
polynomial functions $(F_{\varepsilon} : \mathbb{R} \rightarrow
\mathbb{R})_{\varepsilon}$ such that \ $\mathbb{E} [F_\varepsilon(\varepsilon^{1 / 2}
X_{\varepsilon})] = 0$ then for all $p \geqslant 1$ and under some technical
assumptions we have for every $\psi \in C^\infty_c(\mathbb R,\mathbb R)$
\[ {\psi (t)} F_{\varepsilon} (\varepsilon^{1 / 2} X_{\varepsilon} (t, x))
   \rightarrow 0 \]
almost surely along subsequences in $B^{- \kappa}_{\infty, \infty}
(\mathbb{R} \times \mathbb{T})$ for some $\kappa > 0$. Moreover, if the component of
$F_{\varepsilon} (\varepsilon^{1 / 2} X_{\varepsilon})$ in
the first chaos of the random field $X_\varepsilon$ vanishes, then the above convergence can be improved
at the cost of reducing the space-time regularity:
\[ \varepsilon^{- 1 / 2 - \kappa} \psi (t) F_{\varepsilon} (\varepsilon^{1 /
   2} X_{\varepsilon} (t, x)) \rightarrow 0 \]
as a space--time distribution of parabolic regularity $- 1 / 2 - \kappa$, almost surely
at least along subsequences. See also~\cite{gubinelli_hairer--quastel_2016} where similar results are derived in a specific stationary setting based on the chaos expansion under the stationary measure of $X_\varepsilon$.

%\section{Singular operators}\label{sec:operators}

\section{Anderson Hamiltonian}
\label{sec:operators1}

Paracontrolled distributions and related tools can not only be used to solve
singular SPDEs, they also allow us to construct certain operators that are a
priori ill-defined. Consider for example the parabolic Anderson model
$\partial_t u = \Delta u + u \xi$, that is gPAM with $G (u) = u$. If we
consider the \tmtextit{Anderson Hamiltonian}
\[ \mathcal{H}u = (\Delta + \xi) u, \]
then formally the solution to the parabolic Anderson model is given by $u (t)
= e^{t\mathcal{H}} u_0$, where $(e^{t\mathcal{H}})_{t \geqslant 0}$ is the
semigroup generated by $\mathcal{H}$. So by understanding $\mathcal{H}$ we
should also gain a better understanding of the parabolic equation. In
particular it will be interesting to study the spectrum of $\mathcal{H}$
and the structure of its eigenfunctions in order to learn something about the long time behavior of $u(t)$.

We would like to see $\mathcal{H}$ as an unbounded operator on $L^2
(\mathbb{T}^d)$ and not $\mathcal{C}^{\alpha} (\mathbb{T}^d)$, because $L^2$
is a Hilbert space and the spectral analysis of operators is
much easier on Hilbert spaces than on Banach spaces. For $d = 1$ Fukushima and Nakao~\cite{fukushima_spectra_1977} constructed
$\mathcal{H}$ already in 1977, but the case $d \in\{2,3\}$
was only very recently understood by Allez and Chouk~\cite{allez_continuous_2015}. In
the following we sketch their results for $d = 2$. Of course, there is no
problem to make sense of $\mathcal{H}u$ if $u$ is a smooth function. The
problem is rather that $\mathcal{H}u$ should be in $L^2$, but for $u \in
C^{\infty}$ the product $\xi u$ will not be better behaved than $\xi$ because the
multiplication with a smooth function does not increase the regularity. Since
$\xi$ is a distribution and not an $L^2$ function, $C^{\infty}$ will not be
contained in the domain of $\mathcal{H}$! On the other side if $u$ is too
irregular, the product $\xi u$ may not be defined. So the idea of Allez and
Chouk is to define a domain of paracontrolled functions $u$ for which
$\mathcal{H}u$ takes values in $L^2$. For $\alpha \in \mathbb{R}$ we write
$H^{\alpha} = B^{\alpha}_{2, 2}$ for the $L^2$-Sobolev space with regularity
$\alpha$ and let
\[ \mathcal{D}= \{ u \in H^{1 -} : u^{\sharp} = u - u \para X \in H^{2 -} \},
\]
where $X$ is to be determined. Now we have to deal with Besov spaces other
than $B^{\alpha}_{\infty, \infty}$, but it is easy to see that we have
analogous estimates for the paraproduct and resonant term on $H^{\alpha}$
spaces, more precisely
\begin{gather*}
  \ast \para \ast : H^{\alpha} \times \mathcal{C}^{\beta} \rightarrow H^{\beta
  \wedge (\beta + \alpha)}, \qquad \qquad \alpha, \beta \in \mathbb{R},\\
  \ast \reso \ast : H^{\alpha} \times \mathcal{C}^{\beta} \rightarrow
  \mathcal{C}^{\beta + \alpha}, \qquad \alpha + \beta > 0.
\end{gather*}
The commutator estimate also extends to more general Besov spaces, see~\cite{promel_rough_2016}: For $\alpha \in (0, 1)$ and $\beta, \gamma
\in \mathbb{R}$ with $\beta + \gamma < 0$ and $\alpha + \beta + \gamma > 0$
we have
\[ C\colon H^{\alpha} \times \mathcal{C}^{\beta} \times \mathcal{C}^{\gamma}
   \rightarrow H^{\alpha + \beta + \gamma} . \]
Recall that $\xi \in \mathcal{C}^{- 1 -}$, so for $X \in \mathcal{C}^{1 -}$
and $u \in \mathcal{D}$ we obtain
\begin{align}\label{eq:Anderson-expansion} \nonumber
   \Delta u + u \xi & = \underbrace{(\Delta (u
  \para X) - u \para \Delta X)}_{H^{0 -}} + \underbrace{\Delta
  u^{\sharp}}_{H^{0 -}} + \underbrace{u \para \Delta X}_{H^{- 1 -}} +
  \underbrace{u \para \xi}_{H^{- 1 -}} + \underbrace{u \lpara \xi}_{H^{0 -}} \\
  &\quad + \underbrace{u^{\sharp} \reso \xi}_{H^{1 -}} + \underbrace{C (u, X,
  \xi)}_{H^{1 -}} + \underbrace{u (X \reso \xi)}_{H^{0 -}} .
\end{align}
Choosing $X$ such that $\Delta X = - \xi$ we get $\mathcal{H}u \in H^{0 -}$
for all $u \in \mathcal{D}$, and therefore $\mathcal{H}$ is an unbounded
operator on $H^{0 -}$ with domain $\mathcal{D}$. Moreover, $\mathcal{D}$ is
dense in $H^{0 -}$ and a Hilbert space when equipped with the norm $\| u
\|_{\mathcal{D}} = \| u \|_{H^{1 -}} + \| u^{\sharp} \|_{H^{2 -}}$. But we are
interested in the spectral theory for $\mathcal{H}$ and a generic $u \in
\mathcal{D}$ can never be an eigenfunction because $\mathcal{H}u \in H^{0 -} \not\subset
\mathcal{D}$. To find the eigenfunctions of $\mathcal{H}$ we should identify a
subspace of $\mathcal{D}$ on which $\mathcal{H}$ has better regularity. The
idea of Allez and Chouk is to consider what they call \tmtextit{strongly
paracontrolled distributions}: From the
expansion~(\ref{eq:Anderson-expansion}) we see that if we want $\mathcal{H}u
\in H^{1 -}$, then all terms of regularity $H^{0 -}$ on the right hand side
should cancel up to a remainder of regularity $H^{1 -}$, so we should have
\[ (\Delta (u \para X) - u \para \xi) + \Delta u^{\sharp} + u \lpara \xi + u (X
   \reso \xi) \in H^{1 -}, \]
or in other words
\[ - \Delta u^{\sharp} = (\Delta (u \para X) - u \para \xi) + u \lpara \xi + u
   (X \reso \xi) + u^R, \]
where $u^R \in H^{1 -}$, from where we get\footnote{Strictly speaking it is not possible to invert the Laplace operator
and we have to shift it and consider $(1 - \Delta)^{- 1}$ instead, but for
simplicity we ignore this here.}
\[ u^{\sharp} = (- \Delta)^{- 1} ((\Delta (u \para X) - u \para \xi) + u \lpara
   \xi + u (X \reso \xi)) + u^{\sharp \sharp} = : \Phi (u) + u^{\sharp
   \sharp}, \]
where $\Phi (u)$ is defined through the equation and $u^{\sharp \sharp} \in
H^{3 -}$. Thus, the space of strongly paracontrolled
distributions is
\[ \tmop{dom} (\mathcal{H}) = \{ u \in H^{1 -} : u^{\sharp \sharp} = u - u
   \para X - \Phi (u) \in H^{3 -} \} \subset \mathcal{D}, \]
and for $u \in \tmop{dom} (\mathcal{H})$ we get
\[ \mathcal{H}u = \Delta u + u \xi = (\Phi (u) + u^{\sharp \sharp}) \reso \xi
   + C (u, X, \xi) \in H^{1 -} \subset L^2 . \]
It is not at all trivial that $\tmop{dom} (\mathcal{H})$ contains more
functions than just $0$, but in~\cite{allez_continuous_2015} it is even shown that
$\tmop{dom} (\mathcal{H})$ is dense in $L^2$. Moreover, it is shown that
$\mathcal{H}$ is a symmetric operator:
\[ \langle \mathcal{H}u, v \rangle_{L^2} = \langle u, \mathcal{H}v
   \rangle_{L^2} \]
for all $u, v \in \tmop{dom} (\mathcal{H})$, and there exists a decreasing
sequence $\lambda_1 \geqslant \lambda_2 \geqslant \ldots$ of real eigenvalues with
$\lim_{n \rightarrow \infty} \lambda_n = - \infty$ and an $L^2$-orthonormal basis of
corresponding eigenvectors $(e_n)_{n \in \mathbb{N}}$ such that
$\mathcal{H}e_n = \lambda_n e_n$ for all $n \in \mathbb{N}$ and such that
\[ \mathcal{H}u = \sum_{n = 1}^{\infty} \lambda_n \langle u, e_n
   \rangle_{L^2}, \qquad u \in \tmop{dom} (\mathcal{H}) . \]
There is just one thing that we omitted: as for the parabolic Anderson model it is
of course necessary to renormalise the operator, because if
$(\xi_{\varepsilon})_{\varepsilon > 0}$ is a convolution approximation of
$\xi$ and $- \Delta X_{\varepsilon} = \xi_{\varepsilon}$, then the term
$X_{\varepsilon} \reso \xi_{\varepsilon}$ does not converge but only
$X_{\varepsilon} \reso \xi_{\varepsilon} - c_{\varepsilon}$ converges for a
suitable sequence of diverging constants $(c_{\varepsilon})_{\varepsilon >
0}$. Replacing $X \reso \xi$ by $X \reso \xi - \infty$ has the effect of
changing $\mathcal{H}u = \Delta u + u \xi$ to
\[ \Delta u + u (\xi - \infty) = \Delta u + u \para \xi + u \lpara \xi +
   u^{\sharp} \reso \xi + C (u, X, \xi) + u (X \reso \xi - \infty) . \]
From here we see that if $\xi$ is the space white noise, the operator
$\mathcal{H}$ cannot be continuously extended from $\tmop{dom} (\mathcal{H})$
to the smooth functions because for $u \in C^{\infty}$ the product $u \xi$ does
not create any divergences so $u (\xi - \infty)$ does not make any sense!

Allez and Chouk~\cite{allez_continuous_2015} then proceed to study how the
largest eigenvalue $\lambda_1$ behaves in the white noise case, and they show
that there exist constants $C_1, C_2 > 0$ such that
\begin{equation}
  \label{eq:exp-tails} e^{- C_1 x} \leqslant \mathbb{P} (\lambda_1 \geqslant
  x) \leqslant e^{- C_2 x}
\end{equation}
for $x \rightarrow \infty$. From here we learn at least heuristically that at
large times the solution $u$ to the parabolic Anderson model should not have
any moments, because
\[ u (t) =  e^{t\mathcal{H}} u_0 = \sum_{n = 1}^{\infty} e^{t \lambda_n} \langle u_0, e_n
   \rangle_{L^2}, \]
and omitting the contribution from all eigenvalues except $\lambda_1$ we get
\[ \mathbb{E} [| e^{t \lambda_1} \langle u_0, e_n \rangle_{L^2} |^p]
   =\mathbb{E} [e^{t p \lambda_1} | \langle u_0, e_n \rangle_{L^2} |^p] . \]
By~(\ref{eq:exp-tails}) we have $\mathbb{E} [e^{t p \lambda_1}] = \infty$ as
soon as $t p > C_1$, so in that case we expect that also $\mathbb{E} [| u (t) |^p] = \infty$.
But it remains an open problem how to make this intuitive argumentation
rigorous.

\section{Singular martingale problem}
\label{sec:operators2}
Similar ideas that we developed to study the Anderson Hamiltonian
$\mathcal{H}u = (\Delta + \xi) u$ have also been used by Cannizzaro and Chouk~\cite{cannizzaro_multidimensional_2015}, inspired by \cite{delarue_rough_2016}, to make sense
of certain diffusions with distributional drift: Let\footnote{Besov spaces on $\mathbb{R}$ are defined exactly in the same way as on $\mathbb{T}$ and they have essentially the same properties.} $\xi \in
\mathcal{C}^{\alpha - 1} (\mathbb{R})$ with $\alpha \in (1 / 3, 1 / 2)$ and
consider the SDE $x\colon \mathbb{R}_+ \rightarrow \mathbb{R}$,
\[ \mathd x_t = \xi (x_t) \mathd t + \sqrt{2} \mathd w_t, \]
where $w$ is a Brownian motion. Of course, this equation does not make any
sense because $\xi$ is a distribution and cannot be evaluated in $x_t$. We
still formally write down the martingale problem and call a continuous
stochastic process $x$ a solution if for all suitable $u$ the process
\[ M^u_t = u (x_t) - u (x_0) - \int_0^t \mathcal{G}u (x_s) \mathd s, \qquad t
   \geqslant 0, \]
is a continuous martingale, where $\mathcal{G}u = \xi \partial_x u + \Delta
u$. But just as for the Anderson Hamiltonian the problem is that
$\mathcal{G}u$ is only a distribution whenever $u$ is a smooth function, and
therefore we first need to identify a suitable domain of functions for which
$\mathcal{G}u$ is continuous. We can do this by solving the equation
\begin{equation}
  \label{eq:martingale-resolvent} (\mathcal{G}- \lambda) u = \varphi,
\end{equation}
for $\varphi \in C_b (\mathbb{R})$, the continuous and bounded functions on
$\mathbb{R}$, and $\lambda > 0$. Then we get for the solution $u$ that
$\mathcal{G}u = \varphi + \lambda u$, so provided that $u$ itself is a
continuous function the martingale problem above makes sense. To
solve~(\ref{eq:martingale-resolvent}) we make the paracontrolled Ansatz $u = u' \para X +
u^{\sharp}$ with $u' \in \mathcal{C}^{\alpha} (\mathbb{R})$, $(\Delta -
\lambda) X = \partial_x \xi \in \mathcal{C}^{\alpha} (\mathbb{R})$, and
$u^{\sharp} \in \mathcal{C}^{2 \alpha} (\mathbb{R})$, and we reformulate the
equation as
\[ (\Delta - \lambda) u = \xi \partial_x u + \varphi . \]
From here it is not difficult to see that for $\lambda$ large enough
(depending only on $\xi$ and $X \reso \xi$ but not on $\varphi$) there exists
a unique paracontrolled solution $u$ to the equation, and the space of
paracontrolled functions $u$ which solve $(\mathcal{G}- \lambda) u = \varphi$ for some
$\varphi \in C_b (\mathbb{R})$ is a domain for $\mathcal{G}$. Moreover, there
exists a unique (in law) solution $x$ to the martingale problem defined above.
For details see~\cite{cannizzaro_multidimensional_2015}.

%\bibliographystyle{alpha}
%\bibliography{roeckner-gubinelli-bib.bib}

\begin{thebibliography}{MWX16}

\bibitem[AC15]{allez_continuous_2015}
Romain Allez and Khalil Chouk.
\newblock The continuous {A}nderson {H}amiltonian in dimension two.
\newblock {\em arXiv preprint arXiv:1511.02718}, 2015.

\bibitem[AR91]{albeverio_stochastic_1991}
Sergio Albeverio and Michael R{\"o}ckner.
\newblock Stochastic differential equations in infinite dimensions: solutions
  via {D}irichlet forms.
\newblock {\em Probability theory and related fields}, 89(3):347--386, 1991.

\bibitem[BB16a]{bailleul_heat_2016}
Isma{\"e}l Bailleul and Frederic Bernicot.
\newblock Heat semigroup and singular {PDE}s.
\newblock {\em Journal of Functional Analysis}, 270(9):3344--3452, 2016.

\bibitem[BB16b]{bailleul_higher_2016}
Isma{\"e}l Bailleul and Fr{\'e}d{\'e}ric Bernicot.
\newblock Higher order paracontrolled calculus.
\newblock {\em arXiv preprint arXiv:1609.06966}, 2016.

\bibitem[BCD11]{bahouri_fourier_2011}
Hajer Bahouri, Jean-Yves Chemin, and Raphaël Danchin.
\newblock {\em Fourier {Analysis} and {Nonlinear} {Partial} {Differential}
  {Equations}}.
\newblock Springer, January 2011.

\bibitem[BDH16]{bailleul_quasilinear_2016}
Ismael Bailleul, Arnaud Debussche, and Martina Hofmanova.
\newblock Quasilinear generalized parabolic {A}nderson model.
\newblock {\em arXiv preprint arXiv:1610.06726}, 2016.

\bibitem[BHZ16]{bruned_algebraic_2016}
Yvain Bruned, Martin Hairer, and Lorenzo Zambotti.
\newblock Algebraic renormalisation of regularity structures.
\newblock {\em arXiv:1610.08468 [math]}, October 2016.
\newblock arXiv: 1610.08468.

\bibitem[Bon81]{bony_calcul_1981}
Jean-Michel Bony.
\newblock Calcul symbolique et propagation des singularit{\'e}s pour les
  {\'e}quations aux d{\'e}riv{\'e}es partielles non lin{\'e}aires.
\newblock In {\em Annales scientifiques de l'{\'E}cole Normale sup{\'e}rieure},
  volume~14, pages 209--246, 1981.

\bibitem[CC13]{catellier_paracontrolled_2013}
R{\'e}mi Catellier and Khalil Chouk.
\newblock Paracontrolled distributions and the 3-dimensional stochastic
  quantization equation.
\newblock {\em arXiv preprint arXiv:1310.6869}, 2013.

\bibitem[CC15]{cannizzaro_multidimensional_2015}
Giuseppe Cannizzaro and Khalil Chouk.
\newblock Multidimensional {SDE}s with singular drift and universal
  construction of the polymer measure with white noise potential.
\newblock {\em arXiv preprint arXiv:1501.04751}, 2015.

\bibitem[CF14]{chouk_support_2014}
Khalil Chouk and Peter~K Friz.
\newblock Support theorem for a singular semilinear stochastic partial
  differential equation.
\newblock {\em arXiv preprint arXiv:1409.4250}, 2014.

\bibitem[CFG17]{cannizzaro_malliavin_2017}
Giuseppe Cannizzaro, Peter~K Friz, and Paul Gassiat.
\newblock Malliavin calculus for regularity structures: The case of {gPAM}.
\newblock {\em Journal of Functional Analysis}, 272(1):363--419, 2017.

\bibitem[CGP16]{chouk_invariance_2016}
Khalil Chouk, Jan Gairing, and Nicolas Perkowski.
\newblock An invariance principle for the two-dimensional parabolic {A}nderson
  model with small potential.
\newblock {\em arXiv preprint arXiv:1609.02471}, 2016.

\bibitem[CH16]{chandra_analytic_2016}
Ajay Chandra and Martin Hairer.
\newblock An analytic {BPHZ} theorem for regularity structures.
\newblock {\em arXiv:1612.08138 [math-ph]}, December 2016.
\newblock arXiv: 1612.08138.

\bibitem[DD16]{delarue_rough_2016}
Fran{\c{c}}ois Delarue and Roland Diel.
\newblock Rough paths and 1d {SDE} with a time dependent distributional drift:
  application to polymers.
\newblock {\em Probability Theory and Related Fields}, 165(1-2):1--63, 2016.

\bibitem[DPD03]{daprato_strong_2003}
Giuseppe Da~Prato and Arnauds Debussche.
\newblock Strong solutions to the stochastic quantization equations.
\newblock {\em The Annals of Probability}, 31(4):1900--1916, 2003.

\bibitem[FG16]{furlan_paracontrolled_2016}
M~Furlan and M~Gubinelli.
\newblock Paracontrolled quasilinear {SPDE}s.
\newblock {\em arXiv preprint arXiv:1610.07886}, 2016.

\bibitem[FH14]{friz_course_2014}
Peter~K. Friz and Martin Hairer.
\newblock {\em A {Course} on {Rough} {Paths}: {With} an {Introduction} to
  {Regularity} {Structures}}.
\newblock Springer, August 2014.

\bibitem[FH16]{funaki_coupled_2016}
Tadahisa Funaki and Masato Hoshino.
\newblock A coupled {KPZ} equation, its two types of approximations and
  existence of global solutions.
\newblock {\em arXiv preprint arXiv:1611.00498}, 2016.

\bibitem[FN77]{fukushima_spectra_1977}
Masatoshi Fukushima and Shintaro Nakao.
\newblock On spectra of the {S}chr{\"o}dinger operator with a white {G}aussian
  noise potential.
\newblock {\em Zeitschrift f{\"u}r Wahrscheinlichkeitstheorie und Verwandte
  Gebiete}, 37(3):267--274, 1977.

\bibitem[GIP15]{gubinelli_paracontrolled_2015}
Massimiliano Gubinelli, Peter Imkeller, and Nicolas Perkowski.
\newblock Paracontrolled distributions and singular {PDEs}.
\newblock {\em Forum of Mathematics. Pi}, 3:e6, 75, 2015.

\bibitem[GP15]{gubinelli_lectures_2015}
M.~Gubinelli and N.~Perkowski.
\newblock Lectures on singular stochastic {PDEs}.
\newblock {\em Ensaios Matem\'aticos}, 29, 2015.
\newblock arXiv: 1502.00157.

\bibitem[GP16]{gubinelli_hairer--quastel_2016}
Massimiliano Gubinelli and Nicolas Perkowski.
\newblock The {Hairer}--{Quastel} universality result in equilibrium.
\newblock {\em arXiv:1602.02428 [math-ph]}, February 2016.
\newblock arXiv: 1602.02428.

\bibitem[GP17]{gubinelli_kpz_2017}
Massimiliano Gubinelli and Nicolas Perkowski.
\newblock {KPZ} reloaded.
\newblock {\em Communications in Mathematical Physics}, 349(1):165--269, 2017.

\bibitem[Gub04]{gubinelli_controlling_2004}
M.~Gubinelli.
\newblock Controlling rough paths.
\newblock {\em Journal of Functional Analysis}, 216(1):86--140, 2004.

\bibitem[Hai11]{hairer_rough_2011}
M.~Hairer.
\newblock Rough stochastic {PDEs}.
\newblock {\em Communications on Pure and Applied Mathematics},
  64(11):1547--1585, 2011.

\bibitem[Hai14]{hairer_theory_2014}
M.~Hairer.
\newblock A theory of regularity structures.
\newblock {\em Inventiones mathematicae}, 198(2):269--504, March 2014.

\bibitem[HQ15]{hairer_class_2015}
Martin Hairer and Jeremy Quastel.
\newblock A class of growth models rescaling to {KPZ}.
\newblock {\em arXiv:1512.07845 [math-ph]}, December 2015.
\newblock arXiv: 1512.07845.

\bibitem[Kup16]{kupiainen_renormalization_2016-1}
Antti Kupiainen.
\newblock Renormalization {Group} and {Stochastic} {PDEs}.
\newblock {\em Annales Henri Poincaré. A Journal of Theoretical and
  Mathematical Physics}, 17(3):497--535, 2016.

\bibitem[LCL07]{lyons_differential_2007}
Terry~J. Lyons, Michael~J. Caruana, and Thierry Lévy.
\newblock {\em Differential {Equations} {Driven} by {Rough} {Paths}: {Ecole}
  d'{Eté} de {Probabilités} de {Saint}-{Flour} {XXXIV}-2004}.
\newblock Springer, 1 edition, June 2007.

\bibitem[LQ02]{lyons_system_2002}
Terry Lyons and Zhongmin Qian.
\newblock {\em System {Control} and {Rough} {Paths}}.
\newblock Oxford University Press, 2002.

\bibitem[Lyo98]{lyons_differential_1998}
Terry Lyons.
\newblock Differential equations driven by rough signals.
\newblock {\em Revista Matemática Iberoamericana}, pages 215--310, 1998.

\bibitem[Mey81]{meyer_remarques_1981}
Yves Meyer.
\newblock Remarques sur un théorème de {J}.-{M}. {Bony}.
\newblock In {\em Rendiconti del {Circolo} {Matematico} di {Palermo}. {Serie}
  {II}}, pages 1--20, 1981.

\bibitem[MW15]{mourrat_global_2015}
Jean-Christophe Mourrat and Hendrik Weber.
\newblock Global well-posedness of the dynamic $\phi^4$ model in the plane.
\newblock {\em arXiv preprint arXiv:1501.06191}, 2015.

\bibitem[MW16]{mourrat_global_2016}
Jean-Christophe Mourrat and Hendrik Weber.
\newblock Global well-posedness of the dynamic $\phi^4_3$ model on the torus.
\newblock {\em arXiv preprint arXiv:1601.01234}, 2016.

\bibitem[MWX16]{mourrat_construction_2016}
Jean-Christophe Mourrat, Hendrik Weber, and Weijun Xu.
\newblock Construction of $\phi^4_3$ diagrams for pedestrians.
\newblock {\em arXiv preprint arXiv:1610.08897}, 2016.

\bibitem[OW16]{otto_quasilinear_2016}
Felix Otto and Hendrik Weber.
\newblock Quasilinear {SPDE}s via rough paths.
\newblock {\em arXiv preprint arXiv:1605.09744}, 2016.

\bibitem[PT16]{promel_rough_2016}
David~J Pr{\"o}mel and Mathias Trabs.
\newblock Rough differential equations driven by signals in {B}esov spaces.
\newblock {\em Journal of Differential Equations}, 260(6):5202--5249, 2016.

\bibitem[RZZ15]{rockner_restricted_2015}
Michael R{\"o}ckner, Rongchan Zhu, and Xiangchan Zhu.
\newblock Restricted {M}arkov uniqueness for the stochastic quantization of
  {$P(\Phi)_2$} and its applications.
\newblock {\em arXiv preprint arXiv:1511.08030}, 2015.

\bibitem[ZZ14]{zhu_approximating_2014}
Rongchan Zhu and Xiangchan Zhu.
\newblock Approximating three-dimensional {N}avier-{S}tokes equations driven by
  space-time white noise.
\newblock {\em arXiv preprint arXiv:1409.4864}, 2014.

\bibitem[ZZ15]{zhu_wong_2015}
Rongchan Zhu and Xiangchan Zhu.
\newblock A {W}ong-{Z}akai theorem for $\phi^4_3$ model.
\newblock {\em arXiv preprint arXiv:1504.04143}, 2015.

\end{thebibliography}

\end{document}